\date{} 
\title{Quantitative noise sensitivity\\and exceptional times for percolation}
\author{Oded Schramm \and Jeffrey E. Steif
\thanks{Research supported by the
     Swedish Natural Science Research Council
and the G\"{o}ran Gustafsson Foundation (KVA).}}

\documentclass[12pt,naturalnames]{article}
\usepackage{amsmath}
\usepackage{amssymb}
\usepackage{amsthm}
\usepackage{amsfonts}
\usepackage{graphicx}
\input labelfig.tex
\input epsf.tex

\newif\ifhyper\IfFileExists{hyperref.sty}{\hypertrue}{\hyperfalse}

\ifhyper\usepackage{hyperref}\fi

\newif\ifdraft
\drafttrue

\long\def\comment#1{}
\long\def\oldVersion#1{}

\numberwithin{equation}{section}
\numberwithin{figure}{section}

\newtheorem{theorem}{Theorem}
\numberwithin{theorem}{section}
\newtheorem{corollary}[theorem]{Corollary}
\newtheorem{lemma}[theorem]{Lemma}
\newtheorem{proposition}[theorem]{Proposition}

\theoremstyle{remark}\newtheorem{definition}[theorem]{Definition}
\theoremstyle{remark}\newtheorem{remark}[theorem]{Remark}
\def\eref#1{(\ref{#1})}
\def\QED{\qed\medskip}

\newcommand{\R}{\mathbb{R}}
\newcommand{\C}{\mathbb{C}}
\newcommand{\Z}{\mathbb{Z}}
\newcommand{\N}{\mathbb{N}}
\def\diam{\mathop{\mathrm{diam}}}

\def\dist{\mathop{\mathrm{dist}}}

\def\calC{{\cal C}}
\def\calN{{\cal N}}
\def\bPsi{\mbox{\boldmath $\Psi$}} 

\def\FC(#1,#2){{ 0 \stackrel{#2}{\leftrightarrow} #1 }}
\def\Ito/{It\^o}
\def \eps {\epsilon}
\def \P {{\bf P}}
\def\md{\mid}
\def\Bb#1#2{{\def\md{\bigm| }#1\bigl[#2\bigr]}}
\def\BB#1#2{{\def\md{\Bigm| }#1\Bigl[#2\Bigr]}}
\def\Bs#1#2{{\def\md{\mid}#1[#2]}}
\def\Pb{\Bb\P}
\def\Eb{\Bb\E}
\def\PB{\BB\P}
\def\EB{\BB\E}
\def\Ps{\Bs\P}
\def\Es{\Bs\E}

\def \p {{\partial}}
\def \E {{\bf E}}

\def\var{\operatorname{var}}
\def\closure{\overline}
\def\ev#1{{\mathcal{#1}}}
\def \proof {{ \medbreak \noindent {\bf Proof.} }}
\def\proofof#1{{ \medbreak \noindent {\bf Proof of #1.} }}
\def\proofcont#1{{ \medbreak \noindent {\bf Proof of #1, continued.} }}
\def\CCa{C_3}
\def\CCb{C_1}
\def\CCc{C_2}

\def\bl{\bigl}\def\br{\bigr}

\def\Vss{V_{s,R}\cap V_{s',R}}

\def\nn{[n]}
\def\cro{Q}
\def\croRr{f_r^R}
\def\mfh{\mathfrak H}
\def\hmfh{\hat\mfh}
\def\interior{\operatorname{interior}}
\def\II{{J}}

\def\noopsort#1{}

\begin{document}
\maketitle

 \begin{abstract}
One goal of this paper is to
 prove that dynamical critical site percolation on the planar triangular lattice
has exceptional times at which percolation occurs. 
In doing so, new {\sl quantitative} noise sensitivity results
for percolation are obtained.
The latter is based on a novel method for controlling the
 ``level $k$'' Fourier coefficients via the construction of a
randomized 
algorithm which looks at random bits, outputs the value of a 
particular function but looks at any fixed input bit with low 
probability. We also obtain upper and lower bounds on the Hausdorff
dimension of the set of percolating times. We then study the problem
of exceptional times for certain
``$k$-arm'' events on wedges and cones. As a corollary
of this analysis, we prove, among other things, that there are no times 
at which there are two infinite ``white'' clusters, obtain an upper bound 
on the Hausdorff dimension of the set of times at which there are 
both an infinite white cluster and an infinite black cluster
and prove that for dynamical critical bond percolation on the square grid
there are no exceptional times at which $3$ disjoint infinite clusters are
present.
\end{abstract}

\newpage
\tableofcontents
\newpage

\section {Introduction}

Consider bond percolation on an infinite connected locally finite graph
$G$, where for some $p\in[0,1]$, each edge (bond) of $G$ is, independently of 
all others, open with probability $p$ and closed with probability $1-p$. 
Write $\pi_p$ for this product measure. 
The main questions in percolation theory (see \cite{Grimmett})
deal with the possible existence of infinite connected components 
(clusters) in the random subgraph of $G$
consisting of all sites and all open edges. 
Write $\calC$ for the event
that there exists such an infinite cluster. By Kolmogorov's 0-1 law, 
the probability of $\calC$ is, for fixed $G$ and $p$, either 0 or 1. 
Since $\pi_p(\calC)$ is nondecreasing in $p$, there 
exists a critical probability $p_c=p_c(G)\in[0,1]$ such that
\[
\pi_p(\calC)=\left\{
\begin{array}{ll}
0 & \mbox{for } p<p_c \\
1 & \mbox{for } p>p_c.
\end{array} \right. 
\]
At $p=p_c$ we can have either $\pi_p(\calC)=0$ or $\pi_p(\calC)=1$, depending on $G$.

H\"{a}ggstr\"{o}m, Peres and Steif~\cite{HPS} initiated the study of 
dynamical percolation.
(The notion of dynamical percolation
was invented independently by I.~Benjamini.
While the present paper was motivated by \cite{HPS},
the question studied here had previously been asked by Benjamini, as
we recently became aware.)
In this model, with $p$ fixed, the edges of $G$ switch back 
and forth according to independent 2 state
continuous time Markov chains where closed switches to open at rate $p$ 
and open switches to closed at rate $1-p$. Clearly $\pi_p$ is a stationary
distribution for this Markov process. The general question studied in
\cite{HPS} was whether, when we start with distribution $\pi_p$,
there could exist atypical times at which the percolation 
structure looks markedly different than at a fixed time.

Write $\bPsi_p$ for the underlying probability measure
of this Markov process, and write $\calC_t$ for
the event that there is an infinite cluster of open edges 
at time $t$.

Two results in \cite{HPS} which are relevant to us are

\begin{proposition} \label{pr:noncrit}
 For any graph $G$ we have 
\begin{equation} \left\{ \begin{array}{ccl}
\bPsi_p(\, \calC_t \, \mbox{ occurs for every  } \, t \, )=1  & \mbox{ if } & p>p_c(G) 
   \\[1ex]
\bPsi_p\bigl((\neg\, \calC_t) \mbox{ occurs for every } t\bigr)=1 & \mbox{ if } & p<p_c(G) \, .
\end{array} \right.
\nonumber
\end{equation} 
\end{proposition}

\begin{theorem} \label{th:Zd}
For $d\ge 19,$ the integer lattice $\Z^d$ satisfies
$$
\bPsi_{p_c}\bl((\neg \calC_t) \mbox{ occurs for every } t\br)=1.
$$ 
\end{theorem}

One important aspect of the proof of the latter result is that it uses the fact,
proved in \cite{HS},
that for $d\ge 19$, 
\begin{equation}\label{e.z19}
\pi_p(0\text{ is in an infinite open cluster})
=O(|p-p_c|).
\end{equation}
It is proved  in \cite{KZ} that~\eref{e.z19} does not hold
for $d=2$.
Therefore, the question of whether Theorem 
\ref{th:Zd} is true for $d=2$ becomes interesting. At this point, we mention
that site percolation is the analogous model where the vertices (rather than the edges)
are open or closed independently each with probability $p$ and
dynamical percolation is defined in a completely analogous manner.
Our main result says that Theorem~\ref{th:Zd} does not hold
for site percolation on  the planar triangular grid.
The triangular grid is the graph whose vertex set
is the subset of $\C=\R^2$ consisting of the points
$$
\Z+\exp(2\,\pi\,i/3)\,\Z=
\bigl\{(k+\ell/2,\sqrt 3\, \ell/2):k,\ell\in\Z\bigr\}
$$
and two such points have an edge between them if and only if their distance is 1.
Explicitly stated, our main result is

\begin{theorem} \label{th:except}
Almost surely, the set of times
$t\in[0,1]$ such that dynamical critical site percolation
on the triangular lattice has an infinite open cluster is nonempty.
\end{theorem}

There are no other transitive graphs for which it is known that
dynamical critical percolation has such exceptional times.
(In~\cite{HPS}, it was argued that the event discussed in Theorem~\ref{th:except}
is measurable. A similar comment applies to our other results below.
Thus, measurability issues
will not concern us here.)

We are convinced that Theorem~\ref{th:except}
is true for bond percolation on the square lattice. However, our proof uses 
the existence and exact values of certain so-called 
{\sl critical exponents}, which
are only known to hold for site percolation on the triangular lattice.
These are believed to be the same for 
bond percolation on the square lattice, but even
their existence has not yet been established in that case. 
However, the methods of this paper seem to come quite close
to a proof for the square grid as well: it seems that there are several
ways in which this can perhaps be achieved without determining these
critical exponents. These issues will be further discussed in Section~\ref{sec:square}.

It is interesting to note that by~\cite[Corollary 4.2]{HPS}, a.s.\ at every
time $t$ the set of vertices that are contained in some infinite cluster
has zero density.
\bigskip

On a heuristic level, for Theorem~\ref{th:except} to hold,
it is necessary that
the configuration ``changes fast'' in order to have ``many chances''
to percolate  so that we will in fact have a percolating time. 
Mathematically, ``changing fast'' can be interpreted as having small
correlations over short time intervals, which then suggests the
use of the second moment method which we indeed will use. In other words,
one needs to know that the configuration at a given time tells us almost
nothing about how it will look a short time later. The notion of
``noise sensitivity'' introduced in \cite{BKS} 
is the relevant tool which
describes this phenomenon. We now briefly explain this.

Given an integer $m$, a subset $A$ of $\{0,1\}^m$ and an $\eps>0$, define
$$
N(A,\eps):=
\var\Bigl[\Pb{(Y_1,\ldots,Y_{m})\in A\md X_1,\ldots,X_{m}}\Bigr]
$$
where 
$\{X_i\}_{1\le i\le m}$ are i.i.d.\ with
$\Pb{X_i=1}=1/2=\Pb{X_i=0}$ and conditional on the $\{X_i\}$'s,
$\{Y_i\}_{1\le i\le m}$ are independent with
$Y_i=X_i$ with probability $1-\epsilon$ and
$Y_i=1-X_i$ with probability $\epsilon$. 

\begin{definition} \label{def:noise}
Let $\{n_m\}_{m\ge 1}$ be an increasing sequence in $\N$ going to $\infty$ and 
let $A_m$ be a subset of $\{0,1\}^{n_m}$ for each $m$. We say that
the sequence $\{A_m\}_{m\ge 1}$ is {\bf noise sensitive} if for every $\epsilon >0$,
\begin{equation}\label{e.noisedef}
\lim_{m\to\infty} N(A_m,\eps)=0.
\end{equation}
\end{definition}

This says that for large $m$ knowing the values of
$X_1,\ldots,X_{n_m}$ gives us almost no information concerning whether
$(Y_1,\ldots,Y_{n_m})\in A_m$. This is not the exact definition of noise sensitivity
given in \cite{BKS} 
but is easily shown to be equivalent; see page 14 in that paper. It is also
shown in \cite{BKS} that if (\ref{e.noisedef}) holds for some 
$\epsilon \in(0, 1/2)$, then it holds for all such $\epsilon$
and in addition that $N(A,\eps)$ is decreasing in $\epsilon$
on $[0,1/2]$.

Let $n_m$ be the number of edges in an $(m+1)\times m$ box in $\Z^2$ and let $A_m$ be the event
of a left to right crossing in such a box.
By duality, $\Pb{A_m}=1/2$ for every $m$
(see~\cite{Grimmett}).
In \cite{BKS}, the following result is proved.

\begin{theorem} \label{th:crossing}
The sequence $\{A_m\}_{m\ge 1}$ is noise sensitive.
\end{theorem}

A by-product of the tools needed to prove
 Theorem \ref{th:except} will imply the following more quantitative
version of Theorem~\ref{th:crossing}, which was conjectured in~\cite{BKS}.

\begin{theorem} \label{th:crossingquantsquare}
There exists $\gamma >0$ so that
$$
\lim_{m\to\infty} N(A_m,m^{-\gamma})=0.
$$
\end{theorem}

We have the same result for the triangular lattice but with a better $\gamma$, since
critical exponents are known in this case.

\begin{theorem} \label{th:crossingquant} For critical
site percolation on the triangular lattice,
let $A_m'$ be the event of the existence of a left-right crossing
in a domain $D$ approximating a square of sidelength $m$.
Then for all $\gamma < 1/8$,
$$
\lim_{m\to\infty} N(A_m',m^{-\gamma})=0.
$$
\end{theorem}

In proving our quantitative noise sensitivity results
(Theorems \ref{th:crossingquantsquare} and \ref{th:crossingquant}
as well as those later on necessary for obtaining Theorem \ref{th:except}),
one of two key steps will be Theorem \ref{t.noise}, which gives 
estimates of certain quantities involving
Fourier coefficients of a function based on the properties of an algorithm
calculating the function;
the other key step will be the construction of an appropriate algorithm.
 Precise definitions of undefined terms 
will be given in Section \ref{sec:dilute},
where the connection with noise sensitivity will also be recalled. 

\begin{theorem}\label{t.noise}
Let $n\in\N$ and set $\Omega=\Omega_n:=\{0,1\}^n$.
Let $f:\Omega\to\R$ be a function.
Suppose that there is a randomized algorithm $A$ for
determining the value of $f$ which examines some of the input bits of
$f$ one by one, where the choice of the next bit examined may depend
on the bits examined so far. Let $\II\subseteq \nn:=\{1,2,\dots,n\}$
be the (random) set of bits examined by the algorithm.
Set $\delta=\delta_A:=\sup\bl\{\Ps{i\in \II}:i\in\nn\br\}$.
Then, for every $k=1,2,\dots$,
the Fourier coefficients of $f$ satisfy
\begin{equation}\label{e.noise}
\sum_{S\subseteq\nn,\,|S|=k}\hat f(S)^2 \le \delta\,k\,\|f\|^2,
\end{equation}
where $\|f\|$ denotes the $L^2$ norm of $f$ with respect to the
uniform probability measure on $\Omega$.
\end{theorem}

This result might have some applications to 
theoretical computer science. We will call $\delta_A$
the {\bf revealment} of the algorithm $A$. 
The restriction of $x$ to $J$ (the set of bits examined by the algorithm)
is a {\bf witness} for the function $f$, in the sense that it determines
$f(x)$. As explained in Section \ref{ss.noisetheorem}, Theorem~\ref{t.noise} extends to some
other types of witnesses.

In the case $k=1$, the inequality~\eref{e.noise} cannot be improved
by more than a factor of $O(1/\log n)$:
there is an example showing this with $\delta\le n^{-1/3}\,\log(n)$,
which appears in~\cite[\S 4]{BSW}.
The paper~\cite{BSW} investigates how small the revealment can be for
a balanced boolean function on $\{0,1\}^n$. When the function is
monotone, it is shown that the revealment cannot be much smaller than
$n^{-1/3}$ and in general it cannot be much smaller than
$n^{-1/2}$. Examples are given there which come within logarithmic
factors of meeting these bounds.

We don't know if~\eref{e.noise} is close to being optimal for $k\gg 1$.
One is tempted to speculate that the inequality can be improved
to $\sum_{|S|\le k}\hat f(S)^2\le O(1)\,k\,\delta\,\|f\|^2$.
We do not know any counterexample to this inequality.
However, the AND function $f(x)=\prod_{j=1}^n x_j$
gives an example where
$$
O(1)\,\sum_{|S|= k}\hat f(S)^2 \ge\sqrt k\, \delta\,\|f\|^2
$$
for $k$ satisfying $|k-n/2|=O(n^{1/2})$.
(It is easy to check that the best revealment 
possible for this $f$ is exactly $(2-2^{1-n})/n$.)

\bigskip

Once Theorem \ref{th:except} is established, it is natural to ask: how
large is the set of ``exceptional'' times at which percolation occurs?
In this direction, we have the following result.

\begin{theorem}\label{th.hd}
The Hausdorff dimension of the set of times 
at which dynamical critical site percolation
on the triangular lattice has an infinite cluster is 
an almost sure constant which lies in $[\frac{1}{6},\frac{31}{36}]$.
\end{theorem}

We conjecture that $\frac{31}{36}$ is the correct answer.
In a different direction, once we know that there are exceptional times at which 
percolation occurs, it is natural to ask how many clusters can exist at these
exceptional times. The following provides the answer. 

\begin{theorem} \label{th.no2clusters}
On the triangular lattice, a.s.\ there are no times at which 
dynamical critical site percolation has 2 or more infinite open clusters.
\end{theorem}

For the square grid, we can only prove
\begin{theorem} \label{th.no3clustersZ2}
On $\Z^2$, a.s.\ there are no times at which 
dynamical critical bond percolation has 3 or more infinite open clusters.
\end{theorem}

In some of the figures, we will represent open sites by white hexagons on the
dual grid, and
closed sites by black hexagons. Thus, percolation clusters
correspond to connected components of the union of the white hexagons.
These will also be called white clusters. Likewise, we may also consider
black clusters, which are connected components of black hexagons.

Asking whether 2 infinite white clusters can coexist
at some time is very different from asking whether 2 infinite clusters of 
{\sl different} colors can coexist at some time. We conjecture that there are 
in fact exceptional times 
at which there is both a white and a black infinite cluster and that the 
Hausdorff dimension of such times is $2/3$. We can however prove the following.

\begin{theorem} \label{th.differentclusters}
On the triangular lattice, a.s.\ the
Hausdorff dimension of the set of times at which there is both an infinite white
cluster and an infinite black
cluster is at most $2/3$.
\end{theorem}

We also have the following two results concerning the upper half plane.

\begin{theorem} \label{th.halfplane}
On the triangular lattice intersected with the upper half plane, a.s.\ 
the Hausdorff dimension of the set of times at which there is an
infinite cluster is at most $5/9$.
\end{theorem}

\begin{theorem} \label{th.halfplane2}
On the triangular lattice intersected with the upper half plane, a.s.\ the
set of times at which there is both an infinite white cluster 
and an infinite black cluster
is empty.
\end{theorem}

Theorems~\ref{th.no2clusters}, \ref{th.differentclusters},
\ref{th.halfplane} and \ref{th.halfplane2}
will follow immediately from generalizations presented
in the last part of the paper, which are concerned with studying dynamical percolation
on two other 2 dimensional objects, namely wedges and cones.
For every $\theta \in (0,\infty)$, we let $W_\theta$ denote the 
wedge of angle $\theta$ and $C_\theta$ denote the cone of angle $\theta$. 
For $C_\theta$, we will require that $\theta$ is a multiple of $\pi/3$.
The precise definitions of these will be given in Section \ref{sec:background}.
First, we mention that for all $\theta$, the critical value for site percolation
on $W_\theta$ and on $C_\theta$ is $1/2$, as for site percolation
on the triangular grid and bond percolation on $\Z^2$. 

The following results provide upper and lower bounds on the critical 
angle for which
there are exceptional times for certain $k$-arm type events as well as provide
estimates for the Hausdorff dimension of the set of exception times for a given
angle. In these results, if an upper bound on the Hausdorff dimension
is negative, this means that the set in question is empty. 

We will only do the case where the arms are alternating in color (and hence
for the case of cones, there will be one or an even number of arms).
We do this partially because
it is easier than the general case and because it is all that is needed
in order to make statements concerning the number of infinite clusters.

By a $k$-arm event, we mean an event of the form ``there are $k$ disjoint 
infinite paths having a specified color sequence''; 
for a wedge, the color sequence is well-defined while
for a cone, it is well-defined up to cyclic permutations.

\begin{theorem} \label{th:wedges}
Fix the wedge $W_\theta$ and for integer $k\ge 1$, let 
$A_{W_\theta}^k$ be the event that
there are $k$ infinite disjoint paths
in $W_\theta$ whose colors alternate. Then a.s.\ the Hausdorff 
dimension, $H_{W_\theta}^k$, of the set of exceptional times at which 
$A_{W_\theta}^k$ occurs satisfies
$$
1-\frac{4k(k+1)\pi}{3\theta} \le  
H_{W_\theta}^k
\le
1-\frac{2k(k+1)\pi}{9\theta}.
$$
In particular, for any $k\ge 1$, there are exceptional times
for the event $A_{W_\theta}^k$
for $\theta > \frac{4k(k+1)\pi}{3}$
and there are no exceptional times
for $\theta < \frac{2k(k+1)\pi}{9}$.
\end{theorem}

\begin{theorem} \label{th:cones}
Fix the cone $C_\theta$ with $\theta$ a multiple of $\pi/3$
and let, for $k= 1$ or $k>1$ even, 
$A_{C_\theta}^k$ be the event that
there are $k$ infinite disjoint paths in $C_\theta$
whose colors alternate (if $k>1$). 
Then a.s.\ the Hausdorff 
dimension, $H_{C_\theta}^k$, of the set of exceptional times at 
which $A_{C_\theta}^k$ occurs satisfies
$$
1-\frac{5\pi}{3\theta} \le  
H_{C_\theta}^1
\le
1-\frac{5\pi}{18\theta}
$$
and for $k\ge 2$
$$
1-\frac{4(k^2-1)\pi}{3\theta} \le  
H_{C_\theta}^k
\le
1-\frac{2(k^2-1)\pi}{9\theta}.
$$
In particular, for $k=1$, there are exceptional times
for the event $A_{C_\theta}^1$ for $\theta > \frac{5\pi}{3}$
and there are no exceptional times for $\theta < \frac{5\pi}{18}$,
while for $k\ge 2$, there are exceptional times
for the event $A_{C_\theta}^k$ for $\theta > \frac{4(k^2-1)\pi}{3}$
and there are no exceptional times for $\theta < \frac{2(k^2-1)\pi}{9}$.
\end{theorem}

Theorem \ref{th:cones} is presumably true for other values of $\theta$ 
provided that a proper definition of $C_\theta$ would be given.

\noindent
Remark:
One should note that the upper bounds on the Hausdorff
dimension given in Theorems~\ref{th.hd},\ref{th.differentclusters},\ref{th.halfplane},
\ref{th:wedges} and \ref{th:cones} are all of the form $1-(4/3)\xi$ where $\xi$ is the 
critical exponent for the given event.

\medskip
There is an abstract theory of L\'evy processes
on groups~\cite{Haw84,Eva89}, which gives a
criterion for a L\'evy process (such as $\omega_t$) to hit a set $A$
(such as the set of configurations which contain an infinite component). 
Basically, to show that $A$ is hit, one needs to prove that there exists a
probability measure $\mu$ on $A$ which has $\|\mu\|_*<\infty$
for an appropriate Hilbert norm $\|\cdot\|_*$, based on the Fourier
transform. It seems that we could use this framework in the present
paper, but that would not essentially simplify the core issues we deal with.
Moreover, it seems that our hands-on approach facilitates some generalizations,
which the L\'evy process theory does not cover, which brings us to
our next remark.
\medskip

The fact that the time between flips has an exponential
distribution is not really essential here,
and the results apply in greater generality.
Let $\omega_t(v)$ denote the indicator function for the
event that at time $t$ the site $v$ is white.
  Basically, all of the results 
concerning existence of exceptional times and lower bounds on Hausdorff 
dimension go through (with essentially the same proofs) in the more general setting 
where we assume that

\begin{enumerate}
\item[(i)] The processes $t\mapsto \omega_t(v)$ are independent (possibly with different
distributions depending on $v$) as $v$ runs over all sites. 
\item[(ii)] $\P[\omega_t(v)=1]=1/2$ for all $t$ and $v$. 
\item[(iii)] There is $c>0$ so that
$$
\bigl|\Eb{(-1)^{\omega_t(v)}(-1)^{\omega_s(v)}}\bigr| \le 1-c\, |t-s|
$$
for all $v$ and all $t$ and $s$ satisfying $|t-s|<c$.
\item[(iv)] For each $v$, the process $\omega_t(v)$  has
right continuous paths a.s. 
\end{enumerate}
(Condition (iv) is just a technical condition to insure that the events that we
consider are measurable.)

For results concerning upper bounds on Hausdorff dimension, the proofs go through assuming 
(i), (ii), (iv) and 
\begin{enumerate}
\item[(v)]  There is a $c>0$ so that
$$\E[\text{number of flips of }\omega_t(v) \text{ during }(t_1,t_2)] \le c (t_2-t_1)
$$
for all $v$ and all $t_1,t_2\in\R$ satisfying $t_1 < t_2<t_1+c$.
\end{enumerate}
However, for simplicity, we stick to the original setup.

\medskip

We mention a few other papers where analogous questions to those studied
in \cite{HPS} have been studied for other models. First, the results in
\cite{HPS} were extended and refined in \cite{PS}. Next, analogous questions
for the Boolean model where the points undergo independent Brownian motions was studied
in \cite{BMW}. Analogous questions for the lattice case for certain
interacting particle systems (where updates are not done in an independent fashion)
are studied in \cite{BrSt}. Finally in \cite{BeSc}, it is shown that
there are 
exceptional two dimensionl slices for the Boolean model in four dimensions.

The rest of the paper is organized as follows. 
In Section \ref{sec:dilute}, we will first provide background
on the Fourier-Walsh expansion of a function defined on $\{0,1\}^n$ as well
as connections with noise sensitivity and then continue on 
to give the proof of Theorem \ref{t.noise} as well as a 
generalization to the case where the algorithm is not required to always determine the value 
of the function $f$. (This will accomodate readers who are only interested in
Theorem \ref{t.noise}.) In Section \ref{sec:background}, we will give necessary background 
concerning percolation including a
discussion of critical site percolation on the triangular lattice as well as
a brief discussion of interfaces and critical exponents. 
In Section \ref{sec:percalg}, we will construct two 
algorithms determining certain events involving critical site percolation on the 
triangular lattice and analyze them to obtain upper bounds on the probability that 
a vertex is looked at during the algorithm. (For readers who only want to read
Theorem \ref{t.noise} and see how to apply it, they can just glance through 
Section \ref{sec:background} and then read Section \ref{sec:percalg}.)
Section~\ref{sec:percalg} also
gives a very detailed discussion of interfaces and
completes the proofs of Theorems \ref{th:crossingquantsquare} and 
\ref{th:crossingquant} by applying Theorem \ref{t.noise}. In 
Section \ref{sec:except}, we give the proof of Theorem \ref{th:except}, and
in Section \ref{sec:hd} we give the proof of Theorem \ref{th.hd}.
(Although the upper bound of $31/36$ given in Theorem \ref{th.hd} is a 
special case of Theorem \ref{th:cones}, we choose to give a
different direct proof of this 
without reference to the work done in Section \ref{sec:wedgesupper}.)
In Section \ref{sec:wedgeslower}, we prove the lower bounds 
on the Hausdorff dimension stated in
Theorems \ref{th:wedges} and \ref{th:cones}. 
In Section \ref{sec:wedgesupper}, we prove the upper 
bounds on the Hausdorff dimension
stated in Theorems \ref{th:wedges} and \ref{th:cones}. 
This will be based on a 
general formula which gives an upper bound on the 
Hausdorff dimension of various random sets (or proves they are empty)
in terms of {\sl influences} (Theorem~\ref{th.hdgeneral}).
We conclude the section by showing that
 Theorems \ref{th.no2clusters},
\ref{th.differentclusters}, \ref{th.halfplane} and \ref{th.halfplane2} 
immediately follow from Theorems \ref{th:wedges} and \ref{th:cones}.
After this, in Section \ref{sec:square},
we prove Theorem~\ref{th.no3clustersZ2}
and explain several plausible
ways in which the proof of Theorem \ref{th:except}
might be extended to bond percolation on the square lattice. 
In Section \ref{sec:questions}, we present some open questions.

Finally, the appendix proves some 
results about (non-dynamical) critical percolation
that are needed for Theorems~\ref{th.no2clusters}--\ref{th:cones}.
The main result is that if $r<r'<r''$, then the probability to have
$j$ crossings in a prescribed color sequence between distances $r$ and $r''$
from $0$ is equal to the product of the corresponding probabilities
between radii $r$ and $r'$ and between radii $r'$ and $r''$, times an error
term (depending on $j$) that is bounded away from $0$ and infinity.
Another consequence is that one gets good control on the positions
of the crossings at the inner and outer radii, as was already
demonstrated by Kesten~\cite[Lemma 2]{K2}. The proofs in the appendix
also establish the corresponding statements for critical bond percolation
in $\Z^2$.

\section{Noise sensitivity of algorithmically dilute functions}\label{sec:dilute}

In this section, we give some background and then prove Theorem \ref{t.noise}.

\subsection{Noise sensitivity background}\label{ss.noisebackground}

For a function $f$ from $\Omega=\Omega_n:=\{0,1\}^n$ to $\R$,
the Fourier-Walsh expansion of $f$ is given by
$ f = \sum_{S \subseteq [n]} \hat{f}(S) \chi_S,$ 
where, $\chi_S(T)=(-1)^{|S \cap T|}$ and $\hat{f}(S)= \int f \chi_S$.
Here and in the following, $\int$ refers to integration with respect to
uniform measure and we identify any vector 
$x\in\Omega_n$ with the subset $\{j\in[n]:x_j=1\}$ of
$[n]=\{1,2,\dots,n\}$.
Consequently, $|x|$ denotes the cardinality of that set; that is,
$|x|=\|x\|_1$ for $x\in\Omega_n$. The $\{\chi_S\}_{S\subseteq [n]}$ are an 
orthornormal basis for the $2^n$ dimensional vector space of functions from
$\Omega_n$ to $\R$. In particular,
$$
\|f\|^2=\sum_{S \subseteq [n]} \hat{f}(S)^2.
$$

We now generalize the definition of $N(A,\eps)$ given in the introduction
to any function $f:\Omega\to\R$ by defining
$$
N(f,\eps):=
\var\Bigl[\Eb{f(Y_1,\ldots,Y_{m})|X_1,\ldots,X_{m}}\Bigr].
$$
It is easy to see that (see page 14 in \cite{BKS}) 
\begin{equation}\label{e.page14}
N(f,\eps)= \sum_{\emptyset\neq S \subseteq [n]} \hat{f}(S)^2 (1-2\epsilon)^{2|S|}.
\end{equation}
This explains the importance of the Fourier-Walsh expansion in the study of
noise sensitivity.

\subsection{Proof of Theorem \ref{t.noise}}\label{ss.noisetheorem}

Before giving the proof, we discuss some heuristics.
One may first believe that an estimate such as~\eref{e.noise} would
be valid because when the algorithm terminates, the value of
$f$ is completely determined, and hence perhaps all the nonzero
Fourier coefficients $\hat f(S)\ne 0$, $S\ne\emptyset$, must
satisfy $S\cap \II\ne\emptyset$. However, this is easily shown not 
to be the case. At the $t$-th step of the algorithm,
after $t$ bits have been determined,
we may consider a new function $f_t$, which is $f$ 
with those determined bits substituted. 
If at the $(t+1)$-th step, the $i$-th bit $\omega_i$
of the input $\omega\in\Omega$ is examined, then
in the passage from $f_t$ to $f_{t+1}$, there is a collapsing
of Fourier coefficients: 
$\hat f_{t+1}(S)=\hat f_t(S)+(-1)^{\omega_i} \hat f_t(S\cup\{i\})$
and $\hat f_{t+1}(S\cup\{i\})=0$
for every $S\subseteq\nn\setminus\{i\}$.
Thus, the coefficient $\hat f_{t+1}(S)$ may vanish 
when some bit $i\in S$ is examined by time $t+1$ or when
some $i\notin S$ is chosen at time $t+1$ and it happens that
$\hat f_t(S)+(-1)^{\omega_i} \hat f_t(S\cup\{i\})=0$.
The latter, which we call \lq\lq collapsing from
above\rq\rq, may seem like a highly nongeneric 
situation. However, we cannot rule it out because
we are primarily interested in very non-generic functions,
namely, functions with values in $\{0,1\}$.
The proof below uses a simple decomposition argument to
handle the possibility of
collapsing from above.

In the following, we let $\tilde{\Omega}$ denote the probability
space that includes the randomness in the input bits of $f$ and the randomness
used to run the algorithm and we let $\E$ denote the
corresponding expectation.
Without loss of generality, elements of 
$\tilde{\Omega}$ can be represented as
$\tilde{\omega}=(\omega,\tau)$ where 
$\omega$ are the random bits and $\tau$ represents the
randomness necessary to run the algorithm.

\proof
Fix $k\ge 1$. Let 
$$
g(\omega):=\sum_{|S|=k} \hat f(S)\,\chi_S(\omega)\,,\qquad
\omega\in\Omega.
$$
The left hand side of~\eref{e.noise} is equal to $\|g\|^2$.
Let $\II\subset\nn$ be the random set of all bits examined by the algorithm.
Let $\ev A$ denote the 
minimal $\sigma$-field for which $\II$ is measurable and every
$\omega_i$, $i\in \II$, is measurable; this can be viewed as the
relevant information gathered by the algorithm.
For any function $h:\Omega\to\R$, let
$h_\II:\Omega\to\R$ denote the random function obtained by substituting
the values of the bits in $\II$. 
More precisely, if $\tilde{\omega}=(\omega,\tau)$
and $\omega'\in\Omega$, then $h_\II(\tilde{\omega})(\omega')$ is
$h(\omega'')$ where $\omega''$ is $\omega$ on $\II(\tilde{\omega})$ and is
$\omega'$ on $[n]\backslash \II(\tilde{\omega})$.
In this way, $h_\II$ is a random variable
on $\tilde{\Omega}$ taking values in the set of mappings from
$\Omega$ to $\R$ and it is immediate that this random variable is
$\ev A$-measurable.
When the algorithm terminates, the unexamined bits in $\Omega$
are unbiased and hence 
$\Eb{h\md \ev A}=\int h_\II(=\hat {h_\II}(\emptyset))$
where $\int$ is defined, as usual, to be integration with
respect to uniform measure on $\Omega$.
It follows that $\Es{h}=\Es{\int h_\II}$.

More generally, if $u:\R\to\R$, then
$(u\circ h)_\II=u\circ h_\II$ 
and hence, as above, $ \Es{u(h)}=\EB{\int u(h_\II)}$. 
In particular, for all $h$,
\begin{equation}\label{e.equiv}
\|h\|^2=\Eb{h^2}=\EB{\int h_\II^2}=\Eb{\|h_\II\|^2}.
\end{equation}

Since the algorithm determines  $f$, it is $\ev A$ measurable,
and we have
\begin{equation*}
\|g\|^2=\Es{g\,f}=\EB{\Eb{g\,f\md \ev A}}=
\EB{f\,\Eb{g\md\ev A}}.
\end{equation*}
Since $\Eb{g\md \ev A}=\hat {g_\II}(\emptyset)$,  Cauchy-Schwarz therefore gives
\begin{equation}\label{e.ept}
\|g\|^2 \le \sqrt{ \Es{\hat g_\II(\emptyset)^2}}\,\|f\|\,.
\end{equation}
We may write,
$$
\Es{\hat g_\II(\emptyset)^2}
=
\Eb{\|g_\II\|^2} -
\EB{\sum_{|S|>0} \hat g_\II(S)^2}.
$$
This and~\eref{e.equiv} with $h=g$ imply that
\begin{multline*}
\Es{\hat g_\II(\emptyset)^2}
\le
\|g\|^2
-\EB{\sum_{|S|=k}\hat g_\II(S)^2}
 \\
=
\sum_{S\subseteq\nn}\hat g(S)^2
-\EB{\sum_{|S|=k}\hat g_\II(S)^2}
=
\sum_{|S|=k}\Eb{\hat g(S)^2-\hat g_\II(S)^2}
\,.
\end{multline*}
It is easily seen that for any function $h$,
$h_\II=\sum_{S}\hat h(S)\,(\chi_S)_\II$.
We apply this with $h=g$.
Since $\hat g(S')=0$ if $|S'|> k$,
it follows that for all
$S\subset\nn$ satisfying $|S|=k$
$$
\hat g_\II(S)=
\begin{cases}
\hat g(S),& S\cap \II=\emptyset, \\
0, &S\cap \II\ne\emptyset\,.
\end{cases}
$$
The above estimate for
$\Es{\hat g_\II(\emptyset)^2}$ therefore gives
$$
\Es{\hat g_\II(\emptyset)^2}
\le
\sum_{|S|=k}\hat g(S)^2\,\Pb{S\cap \II\ne\emptyset}
\le \sum_{|S|=k}\hat g(S)^2\,\sum_{i\in S}\Ps{i\in \II}
\le \|g\|^2\,k\,\delta
\,.
$$
Substituting this estimate in~\eref{e.ept} and squaring the resulting
inequality completes the proof of the theorem.
\QED

The theorem may be easily generalized to situations where the
algorithm does not always determine the value of $f$ precisely;
that is, $f_\II(x)$ still depends on $x\in\Omega$.

Set 
$$
\var_\Omega(f_\II):=\int (f_\II)^2 - \Bigl(\int f_\II\Bigr)^2
=\sum_{S\ne\emptyset}\hat f_\II(S)^2\,,
$$
where the integrations are with respect to the uniform
probability measure on $\Omega$.
Note that $\Eb{\var_\Omega(f_\II)}$ is an indicator for how precisely
the algorithm can be used to approximate $f$;
when $\var_\Omega(f_\II)$ is small, with high conditional probability,
$|\hat f_\II(\emptyset)-f|$ is not too large.

When $\var_\Omega(f_\II)\ne0$, we have to replace the calculation
in the proof of Theorem~\ref{t.noise}
by the following
\begin{multline}
\|g\|^2=
\Eb{g_\II\,f_\II}
=\Eb{ g_\II\,\hat f_\II(\emptyset)}
+ 
\Eb{g_\II\,(f_\II-\hat f_\II(\emptyset))}
\\
\le
\Bigl(\Es{\hat g_\II(\emptyset)^2}\,\Es{\hat f_\II(\emptyset)^2}\Bigr)^{1/2}
+
\Bigl(\Es{ g_\II^2}\,\Es{\var_\Omega( f_\II)}\Bigr)^{1/2}. 
\end{multline}
Since $\Eb{\hat f_\II(\emptyset)^2+\var_\Omega(f_\II)}=\Es{f^2}$,
we have $\Es{\hat f_\II(\emptyset)^2}\le \|f\|^2$.
Thus, the above gives
$$
\|g\|^2 \le
\sqrt{\Eb{\hat g_\II(\emptyset)^2}}\,\|f\|
+
\|g\|\,\sqrt{\Es{\var_\Omega(f_\II)}}\,.
$$
Using the same estimate for $\Eb{\hat g_\II(\emptyset)^2}$ as in the proof
of Theorem~\ref{t.noise}, we  obtain
$$
\|g\|^2
\le
\|g\|\,\|f\|\sqrt{k\,\delta}+
\|g\|\,\sqrt{\Es{\var_\Omega(f_\II)}}\,.
$$
Consequently, squaring both sides gives the following
generalization of~\eref{e.noise}:
\begin{multline}\label{e.gnoi}
\sum_{|S|=k}\hat f(S)^2
\le
\Bigl(\|f\|\sqrt{k\,\delta}+
\sqrt{\Es{\var_\Omega(f_\II)}}\Bigr)^2
\\
\le
2\,k\,\delta\,\|f\|^2+
2\,
\Es{\var_\Omega(f_\II)}\, .
\end{multline}
\bigskip

Theorem~\ref{t.noise} holds more generally than stated.
If $W$ is a random subset of $\nn$,
we say that $W$ is a {\bf witness} for $f:\Omega\to\R$ if
the value of $f$ is determined by its restriction to
$W$.  We say that $W$ is {\bf $\delta$-dilute} if
$\max_{i\in\nn}\Pb{i\in W(\omega)}\le\delta$.
The related notion of short witnesses is of central importance
in Talagrand's epic isoperimetric saga~\cite{Talagrand}.
The proof of Theorem~\ref{t.noise} holds in the more general setting where the
random
set $J$ is a witness
with the property that for all subsets $A\subseteq \nn$,
conditioned on $J=A$ (assuming this has positive probability)
and conditioned on $\omega$ restricted to $A$, the
$\{\omega_i\}_{i\not \in A}$ are uniform i.i.d.\ bits.
As pointed out to us by Asaf Nachmias, 
Theorem~\ref{t.noise} does not hold for arbitrary witnesses,
even if we allow for a multiplicative constant in the
right hand side of~\eref{e.noise}: if you take \lq\lq Recursive Ternary Majority\rq\rq\ on
$n=3^h$ bits, there is a (symmetric) witness having only $2^h$ elements,
yielding a $\delta$ which is $(2/3)^h$; however, the sum of the squares
of the level 1 Fourier coefficients is $(3/4)^h$.

\section{Percolation background and notations}\label{sec:background}

Duality plays a central role in the theory of percolation in two dimensions.
A dual-open path on the triangular grid is defined as a path
in the grid whose vertices are all closed. For the square grid,
a dual-open path is defined as a path in the dual of the square grid that
does not intersect any open edge in the primal grid.
The basic observation is that for site percolation on the 
triangular grid at $p=1/2$ the distribution of the collection of
dual-open paths is the same as the distribution of the collection of primal open paths.
(Sometimes, we use the term ``primal open path'', for an open path, to make
the distinction with the dual-open path clearer.)
For critical bond percolation on the square grid at $p=1/2$, the distribution
of the dual-open paths is the image of the distribution of the
open paths under translation by $(1/2,1/2)$.
This simple duality is one of the important ingredients in the proof by
Kesten that $p_c=1/2$ for these two percolation models~\cite[pg.~53]{Kesten}
and the earlier proof by Harris (see \cite{Harris})
that there is a.s.\ no infinite cluster at $p=1/2$.

For $0\le r <R<\infty$, let $A(r,R)$ denote the event that there is an open
crossing of the annulus $r\le|z|\le R$, namely, an open path connecting a vertex
inside the disk $|z|\le r$ to a vertex in $|z|\ge R$.
Let $\alpha(r,R)$ denote the probability of
${A(r,R)}$ at percolation parameter $p=p_c=1/2$. 
Abbreviate $\alpha(0,R)$ by $\alpha(R)$.
For convenience, we adopt the convention $\alpha(r,R)=1$ whenever $r\ge R$.
The function $\alpha(r,R)$ is essentially multiplicative, in the following
sense: there is a constant $C>0$ such that for every $0\le r_1\le r_2\le r_3<\infty$,
\begin{equation}\label{e.multip}
C^{-1}\,\alpha(r_1,r_3)\le \alpha(r_1,r_2)\,\alpha(r_2,r_3)\le C\,\alpha(r_1,r_3)\,.
\end{equation}
In fact, this holds for critical bond percolation on the square grid as well
as for critical site percolation on the triangular grid.
The (standard) proof of~\eref{e.multip}
is based on the Harris-FKG inequality
and the celebrated Russo-Seymour-Welsh (RSW) theorem
(see~\cite{Grimmett,Kesten}).
A proof of a generalization of~\eref{e.multip} is given in the appendix.
Another consequence of RSW that we will use is the existence of a constant
$c>0$ such that for every $r>0$,
\begin{equation}\label{e.rsw}
c\le \alpha(r,2\,r)\,.
\end{equation}

\begin{figure}
\centerline{\epsfysize=2.5in\epsfbox{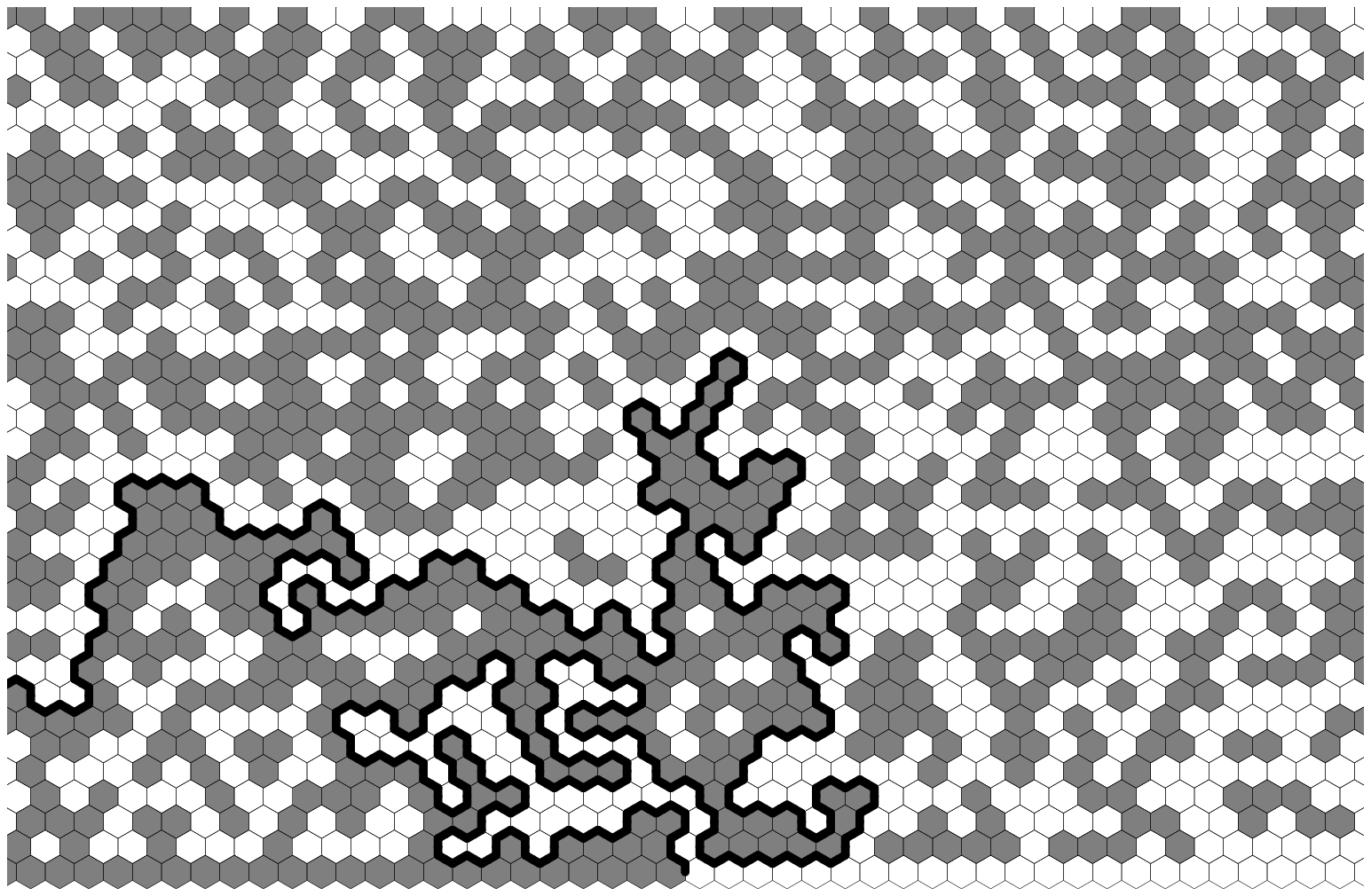}}
\begin{caption} {\label{f.pcurve}The percolation interface.}
\end{caption}
\end{figure}

The Stochastic L\"owner evolution (SLE) introduced in~\cite{Schramm} is a one parameter 
family of random curves indexed by a real positive parameter $\kappa$.
It was conjectured in \cite{Schramm} 
that the scaling limit of outer boundaries of critical
site percolation clusters on the triangular grid
as well as
bond percolation clusters on $\Z^2$ 
are (chordal) SLE$_\kappa$ with ${\kappa=6}$.
Smirnov~\cite{Smirnov,Smirnov1} proved the corresponding statement for critical site
percolation on the triangular lattice. 
(See also \cite{CN}.) We now explain some of this more precisely. We first
perform independent site percolation on the upper half of the triangular lattice 
but declare the sites $\{(k,0):k > 0\}$ to all be open and $\{(k,0):k\le 0\}$ to all be closed.
In the hexagonal grid dual to the triangular grid, there will then be a
unique path in the upper half plane from $(\frac{1}{2},0)$ to $\infty$ 
which has white hexagons containing open 
sites on the right and black hexagons containing closed
sites on the left. (See Figure~\ref{f.pcurve}.)
Smirnov's result is that the limit (in an appropriate toplogy) as the mesh size of the
lattice goes to $0$ of the law of this path
is chordal SLE$_6$.
This path described above, which has open sites
on its right and closed sites on its left, is an example of what is called an 
{\it interface}.

The conformal invariance and the SLE description of critical percolation on 
the triangular lattice allowed researchers to prove a number of conjectures by 
physicists concerning so-called {\sl critical exponents}.
For example, for critical site percolation on the triangular grid, it was established in 
\cite{LSW} that 
\begin{equation}\label{e.onearm}
\alpha(R)=R^{-5/48 + o(1)} \text{ as $R\to\infty$}\,.
\end{equation}
In fact, the same proof actually gives for $R\ge r\ge 1$,
\begin{equation}\label{e.onearmstrong}
\alpha(r,R)=(R/r)^{-5/48 + o(1)} \text{ as $R/r\to\infty$}\,.
\end{equation}

For $1\le r \le R$, the two arm function $\alpha_2(r,R)$ denotes the probability
that there is both an open path from $|z|\le r$ to $|z|\ge R$ and also a 
dual-open path from $|z|\le r$ to $|z|\ge R$. 
We abbreviate $\alpha_2(1,R)$ by $\alpha_2(R)$.

Next, let $M$ be a half plane in $\R^2$ and let $v$ be some vertex in $M$
satisfying $\dist (v,\p M)\le 2$.
Denote by $\alpha^+(R)$ the probability that there is an open
path in $M$ from $v$ to distance at least $R$ away from $v$.
This quantity depends on $R$, $v$, and $M$, but the dependence on
$v$ and $M$ will usually be suppressed.
More generally, for $0\le r< R$, let 
$\alpha^+(r,R)$ denote the probability that there is
an open path in $M$ from some vertex $u$ satisfying $|u-v|\le r$ to some 
vertex $w$ satisfying $|w-v|\ge R$.

As with $\alpha(r,R)$, we adopt the convention
$\alpha^+(r,R)=1$ whenever $r\ge R$.

It is known that
the functions $\alpha^+$ and $\alpha_2$ also satisfy~\eref{e.multip}, with possibly
different constants. (When considering $\alpha^+$,
this applies to any fixed choice of $M$ and $v$.)
These inequalities are valid for site percolation on the triangular grid as well
as bond percolation on the square grid. 
A proof can be found in the appendix.
For site percolation on the triangular grid, the corresponding exponents
were established in~\cite{SW}:
\begin{equation}\label{e.expo}
\alpha^+(R)=R^{-1/3+o(1)},\qquad
\alpha_2(R)=R^{-1/4+o(1)},
\end{equation}
as $R\to\infty$.
In fact, the same proofs actually give for $R\ge r\ge 1$
\begin{equation}\label{e.expostrong}
\alpha^+(r,R)=(R/r)^{-1/3+o(1)},\qquad
\alpha_2(r,R)=(R/r)^{-1/4+o(1)},
\end{equation}
as $R/r\to\infty$.

For bond percolation on the square grid, such exact estimates are unavailable,
because there is currently no proof that the interface converges
to SLE$_6$.
In the case of the square grid, the
estimate $\alpha(r,R)\le C\,(R/r)^{-\eps}$,
where $C,\eps>0$ are constants, follows easily from the 
RSW theorem (see \cite{Grimmett,Kesten}).
The RSW proofs can give an actual value for $\eps$, but it is rather
small. We can also use the obvious estimate $\alpha^+(r,R)\le\alpha(r,R)$
to obtain a similar bound for $\alpha^+$.

\smallskip

We now give the precise definitions for wedges and cones.
For this purpose, we first recall the definition of the infinitely
branched cover of $\R^0$ over $0$.
Let $X=\bl\{(z,\theta):z\in\C\setminus \{0\},\,\theta\in\R,e^{i\theta}|z|=z\br\}$,
and set $\psi(z,\theta)=z$. On the surface $X$ we define the
metric $d_X$ as the pullback of the Euclidean metric of $\R^2$ under
$\psi$, namely, $d_X(x,y)$ is the infimum of the length
of $\psi\circ \gamma$ for any continuous path $\gamma\subset X$ connecting $x$
and $y$. Let $C_\infty$ denote the completion of $(X,d_X)$.
Since $\R^2$ is complete, it is easy to see that
$C_\infty\setminus X$ consists of a single point, which we denote by $0$.
We extend the map $\psi$ by setting $\psi(0)=0$.
Let $V$ be the set of points in $C_\infty$ that are mapped to vertices
of the triangular grid under $\psi$. The triangular grid
on $C_\infty$ has vertices $V$ and an edge between any two vertices
at distance $1$ apart.
Now the wedge $W_\theta\subset C_\infty$ is defined by
$W_\theta:=\{0\}\cup \bl\{(z,\theta')\in X: \theta'\in[0,\theta)\br\}$.
The triangular grid on $W_\theta$ is just the intersection of
the triangular grid on $C_\infty$ with $W_\theta$.

On $C_\infty$ we may define the rotation $R_\theta$
by $R_\theta(0)=0$ and $R_\theta(z,\theta')=(e^{i\theta}z,\theta+\theta')$.
This is clearly an isometry of $C_\infty$.
The cone $C_\theta$ is defined as the quotient $C_\infty/R_\theta$; that is,
the set of equivalence classes of points in $C_\infty$, where two points
are considered equivalent if one is mapped to the other by a power of
$R_\theta$.  Now suppose that $\theta=n\,\pi/3$ where $n\in\N_+$. Then
$R_\theta$ restricts to an isomorphism of the triangular grid on $C_\infty$.
In this case we define the triangular grid on $C_\theta$ as the quotient
of the grid on $C_\infty$ under $R_\theta$.
In other words, the vertices are equivalence classes of vertices in
$C_\infty$ and an edge appears between two equivalence classes if
there is an edge connecting representatives of these classes.
Note that $C_{2\pi}$ is just the euclidean plane with the triangular
lattice.

\smallskip

We end this section with describing the so-called full and half plane 
exponents for $k$-arm events that were derived in \cite{SW}.

For integer $k\ge 1$, let $A^k(r,R)$ be the event that there are $k$ disjoint
crossings of the annulus $\{z\in\R^2: r\le |x|\le R\}$ with a specified color
sequence (up to rotations), where we require that both colors appear 
in the color sequence.
For $k\ge 2$, and $r\ge 10k$, it was proved in \cite{SW} that
\begin{equation}\label{e.fullpk}
\alpha_k(r,R):=\Pb{A^k(r,R)}=(R/r)^{\frac{1-k^2}{12} + o(1)},
\end{equation}
as $R\to\infty$ while $r$ is fixed.
(The result for $\alpha_2(R)$ in~\eref{e.expo} above is a special case of this.)
Next, letting $A_+^k(r,R)$ be the event that there are $k$ disjoint paths in the 
upper half plane from $|z|\le r$ to $|z|\ge R$ with any specified color sequence, then for
$k\ge 1$, and $r\ge 10k$, it was proved in \cite{SW} that
\begin{equation}\label{e.halfpk}
\alpha^+_k(r,R):=\Pb{A_+^k(r,R)}=(R/r)^{\frac{-k(k+1)}{6} + o(1)},
\end{equation}
as $R\to\infty$ while $r$ is fixed. %
(The result for $\alpha_+(R)$ in~\eref{e.expo} is a special case of this.) 
Just as we said that the proofs of~\eref{e.expo} actually yield~\eref{e.expostrong},
it is also the case that the proofs of~\ref{e.fullpk} and~\ref{e.halfpk} also
yield versions when $R/r\to\infty$ while $r\ge 10\,k$ is not necessarily fixed.

\section{Noise sensitivity for percolation} \label{sec:percalg}

\subsection{Simply connected case}\label{ss.simply}

To apply Theorem~\ref{t.noise} to percolation, we will need to describe algorithms achieving
small revealment. One result of that nature is

\begin{theorem}\label{t.square}
Let $\cro=\cro_R$ be the indicator function for the 
event that critical site percolation on the standard triangular grid
contains a left to right crossing in some grid-approximating domain $D$
to a large square of side length $R$.
(For example, we could take $D$ to be the union
of the hexagons in the dual grid that are contained in the
square.)
Then there is a randomized algorithm $A$ determining $\cro$ 
such that $\delta_A \le R^{-1/4+o(1)}$
as $R\to\infty$.

For critical bond percolation on the square grid, there is such an
algorithm satisfying $\delta\le C\,R^{-a}$ for some constants
$a,C>0$.
\end{theorem}

\noindent
Remark:
Theorem~\ref{t.square} says that there is an algorithm for the
relevant event which exposes on average at most $R^{7/4+o(1)}$ bits.
Since the probability of points not too close to the boundary
being pivotal is about $R^{-5/4+o(1)}$ (this is the 4 arm event)
and for a monotone function $f$, $\hat f(\{i\})$ is
the probability that $x_i$ is pivotal,
the case $k=1$ in Theorem~\ref{t.noise} implies that the revealment
is at least $R^{-1/2+o(1)}$.
As pointed out in Peres, Schramm, Sheffield and Wilson~\cite{PSSW}, 
this can also be obtained using an inequality of O'Donnell and Servedio.

Theorems~\ref{t.noise} and \ref{t.square} immediately give

\begin{corollary}\label{c.square}
For every $\eps>0$ there is a constant $C=C(\eps)$ such that
$$
\sum_{|S|=k}\hat\cro_R(S)^2 \le C\, k\,R^{-1/4+\eps}
$$
holds for every $k=1,2,\dots$ and for every $R>0$.
\QED
\end{corollary}

\begin{figure}
\centerline{\epsfysize=2in\epsfbox{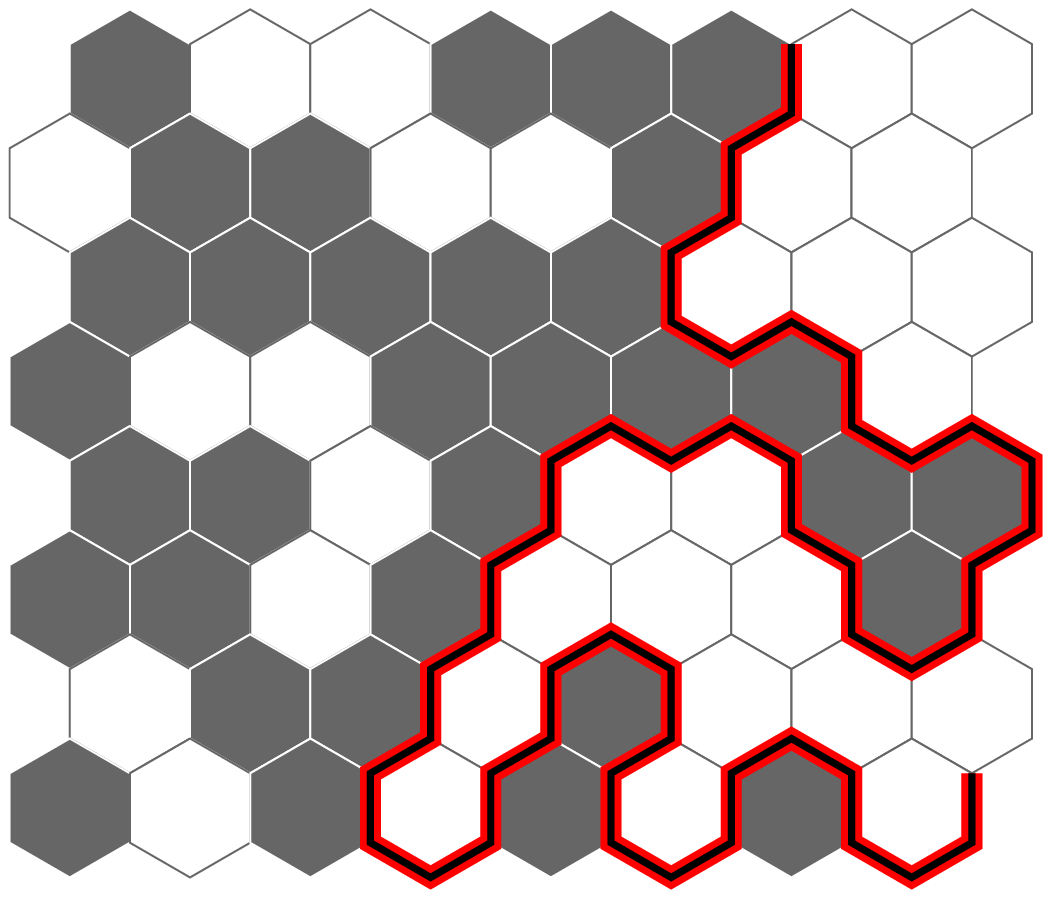}}
\begin{caption} {\label{f.interface}Following the interface from the corner.}
\end{caption}
\end{figure}

The basic idea of the proof of Theorem~\ref{t.square} is rather simple.
First we consider the interface started at the lower right corner and
stopped when it hits the upper or left edges. (See Figure~\ref{f.interface}.)
This interface is sufficient to determine $Q$.
If we traverse the interface, revealing just the bits necessary for its
determination, then with high probability most bits will not be
examined. However, this does not yield an algorithm with small revealment
because the hexagons near the lower right corner are very likely
to be examined. 
To rectify this problem, we instead start the interface at a different
(random) location $p_0$ on the right side
of $D$. 
This determines the existence of a crossing from the right side above
$p_0$ to the left side. Then another interface started at $p_0$ will
determine if there is a crossing that starts below $p_0$.

Let us now be a bit more precise regarding the notion of an interface.
In the following, we use an equivalent dual version of the
site percolation model on the triangular grid. The dual graph is
the hexagonal grid, and we color the hexagon white if the site
contained in it is open and black if the site in it is closed.
Of course, there is no essential difference between these two
representations. The advantage of this dual framework is that the
figures are clearer and the notion of the interface is slightly more
natural.

\begin{figure}
\centerline{\epsfysize=2.5in\epsfbox{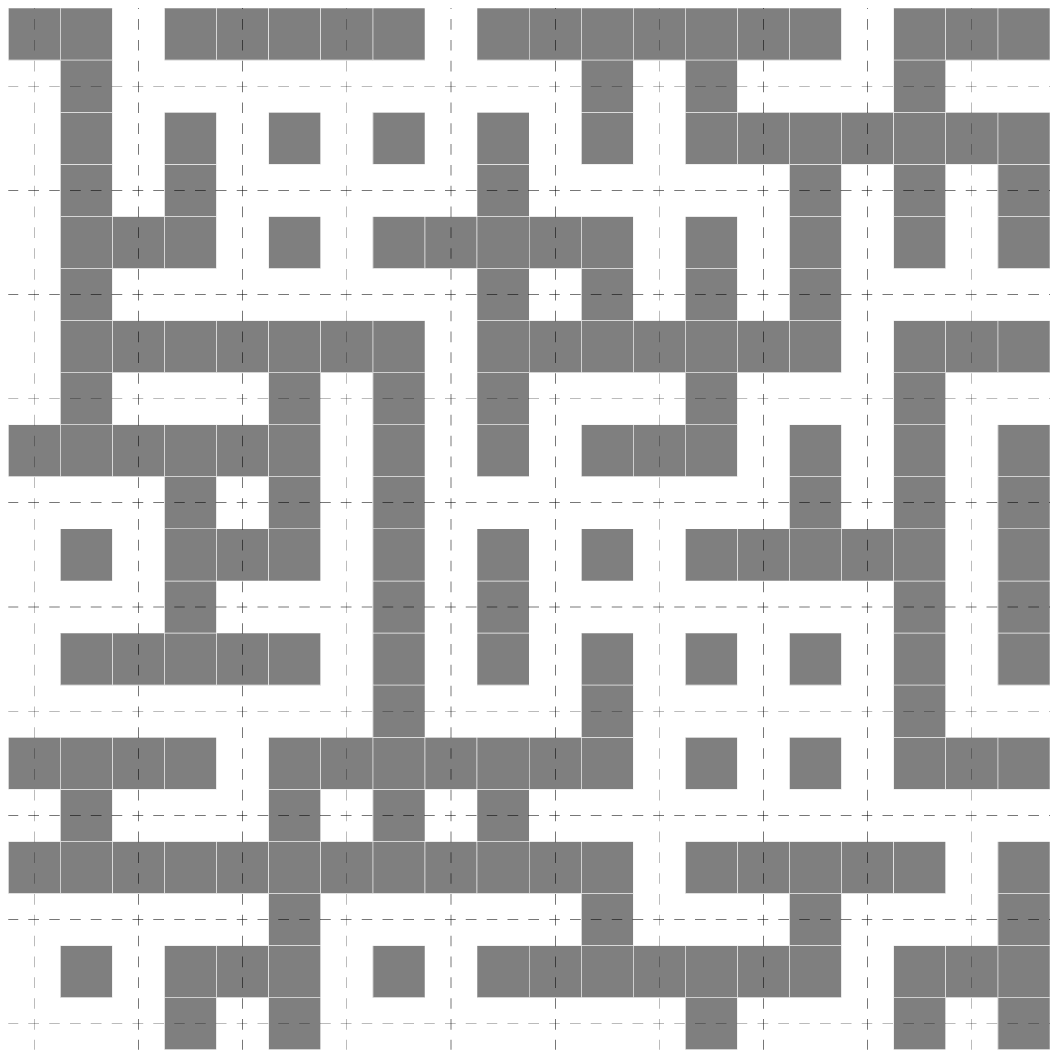}}
\begin{caption} {\label{f.squareColorPerc}A color scheme for bond percolation on $\Z^2$.}
\end{caption}
\end{figure}

Note that bond percolation on the square grid also has a similar
coloring representation. One such scheme is to
color the squares of sidelength $1/2$ centered at the sites
of $\Z^2$ white, color the squares of sidelength $1/2$ that are
concentric with the square faces of $\Z^2$ black,
and color each square of sidelength $1/2$ whose center is the
midpoint of an edge of $\Z^2$ white or black if that edge is
open or closed, respectively. See Figure~\ref{f.squareColorPerc}.
This scheme has the important property that the boundary between
the union of the white clusters is a $1$-manifold; that is,
at every vertex of the grid $(1/2)\,\Z^2+(1/4,1/4)$ there are
two edges that are on the common boundary.

\begin{figure}
\SetLabels
\B(.56*.98)$\zeta$\\
\R(.19*.48)$p_0$\\
\endSetLabels
\centerline{\epsfysize=2in%
\AffixLabels{%
\epsfbox{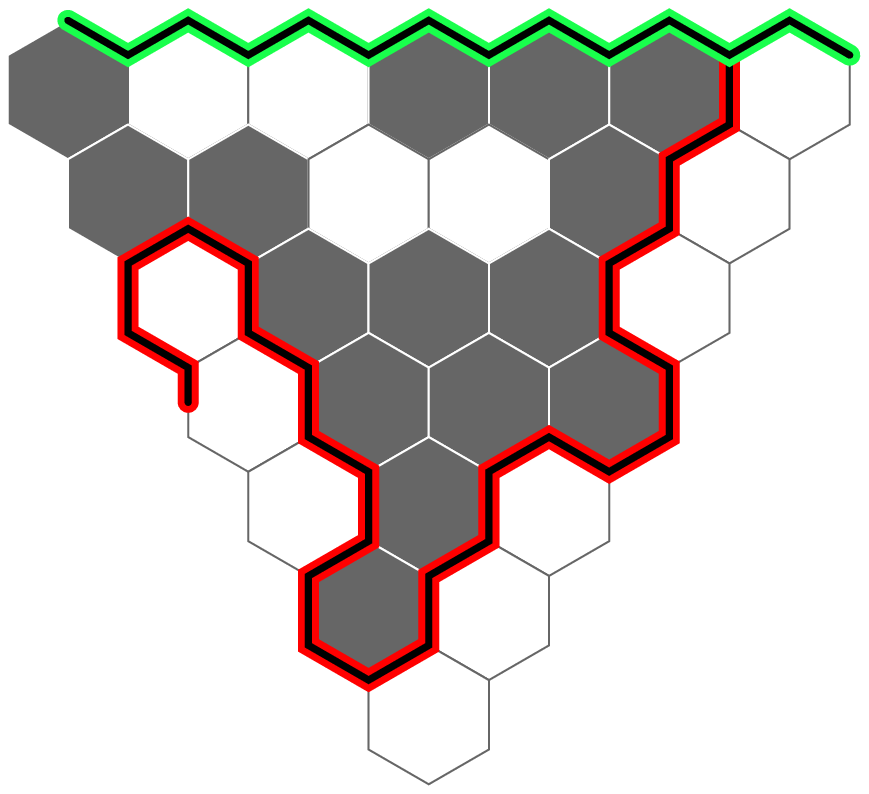}}%
}
\begin{caption} {\label{f.interfaceDef}An interface started at $p_0$
and headed towards $\zeta$.}
\end{caption}
\end{figure}

We now assume that $D$ is a bounded simply connected domain which is
the interior of a union of hexagons in the hexagonal grid.
Suppose that $p_0$ is a point on $\p D$ that is on the boundary
of a single hexagon in $D$, and that $\zeta$ is a closed arc
in $\p D\setminus\{p_0\}$.
The interface in $\closure D$ from $p_0$ to $\zeta$ is a 
random path $\beta$ contained in the $1$-skeleton of the hexagonal grid
starting at $p_0$ and ending at a point in
$\zeta$, as indicated in Figure~\ref{f.interfaceDef}.
We can precisely define $\beta$ as the unique oriented simple path
from $p_0$ to $\zeta$ that is contained in
the union of the boundaries of the hexagons
contained in $\closure D$ and satisfies
(1) $\beta\cap\zeta$ consists of
the terminal point of $\beta$,
(2) whenever $\beta$ traverses an arc along the boundary
of a black hexagon $H\subset\closure D$, the 
arc is traversed counterclockwise around $\p H$, and
(3) whenever $\beta$ traverses an arc along the boundary
of a white hexagon $H\subset\closure D$, the 
arc is traversed clockwise around $\p H$.
It is easy to verify that this uniquely defines $\beta$,
as follows.
First, the initial arc of $\beta$ is determined by the
color of the hexagon in $\closure D$ containing $p_0$.
When $\beta$ first meets a hexagon contained in 
$\closure D$, its turn is clearly specified.
(If the hexagon is black, then $\beta$ must make
a $\pi/3$ turn to the right,
and if the hexagon is white, then $\beta$
must make a $\pi/3$ turn to the left.)
Now consider the situation in which $\beta$
first meets a hexagon that is not contained in
$\closure D$. Let $\beta'$ be the
arc of $\beta$ from $p_0$ up to that point.
Then $\beta$ at that point makes the turn into
the component of $\closure D\setminus \beta'$
that contains $\zeta$, as in the figure.

Another way to describe this interface is that we color 
the counterclockwise arc of $\p D$  from $p_0$ to
$\zeta$ white and the clockwise arc from $p_0$ to
$\zeta$ black, and $\beta$ then is the common boundary
component between white and black in $\closure D$
starting at $p_0$.
We will call this type of interface a {\bf chordal interface},
to differentiate it from the interface that we will later need
when discussing the annulus crossing event.
(The chordal interface was proved by Smirnov~\cite{Smirnov1} to converge to chordal
SLE(6).)

\proofof{Theorem~\ref{t.square}}
We mostly concentrate on the case of site percolation on the
triangular grid. The details in the case of bond percolation
on $\Z^2$ are essentially the same.

\begin{figure}
\centerline{\epsfysize=2in%
\SetLabels
\L(1*.43)$\beta$\\
\L(.95*.6)$p_0$\\
\endSetLabels
\AffixLabels{\epsfbox{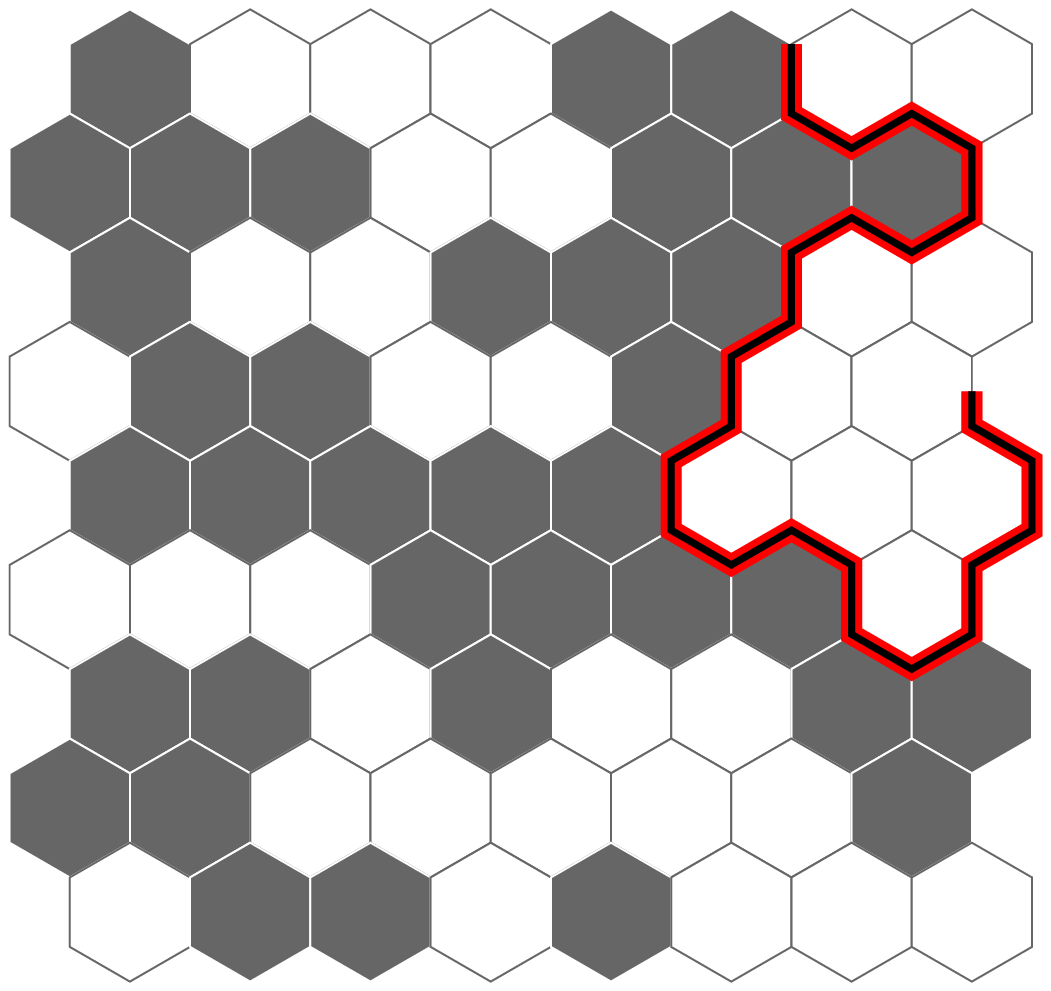}}%
\qquad
\SetLabels
\L(1.01*.25)$\beta'$\\
\L(.95*.6)$p_0$\\
\endSetLabels
\epsfysize=2in%
\AffixLabels{\epsfbox{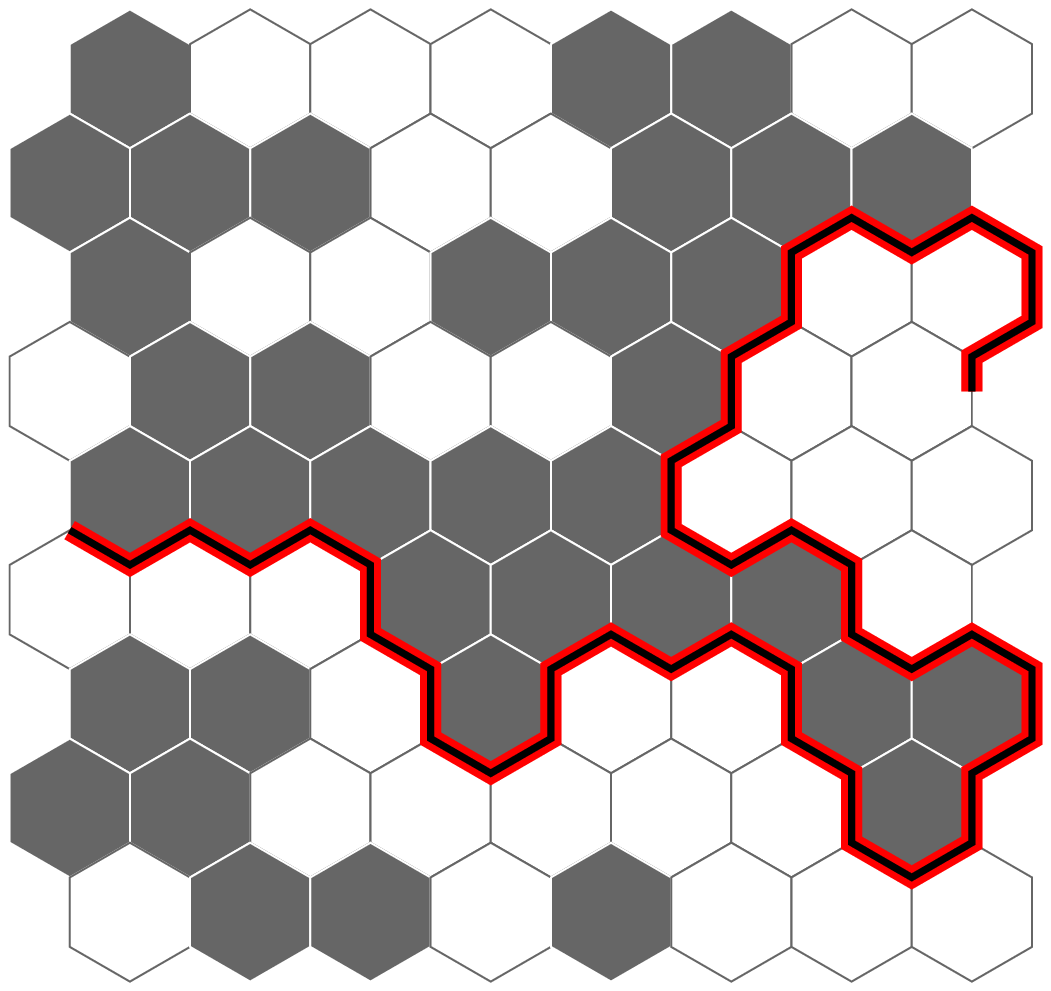}}%
}
\begin{caption} {\label{f.beta}The interfaces $\beta$ and $\beta'$.}
\end{caption}
\end{figure}

The algorithm proceeds as follows. (See Figure~\ref{f.beta}.)
There are four distinguished boundary arcs of $D$,
which we call \lq\lq left\rq\rq, \lq\lq right\rq\rq, \lq\lq up\rq\rq\
and \lq\lq down\rq\rq.
Pick uniformly at random 
an edge $e_0$ on the right hand boundary of 
$D$, and let $p_0$ be its midpoint.
Let $\zeta$ be the union of the top and left boundary
segments of $D$.
Explore the interface $\beta$ from $p_0$ to $\zeta$, examining the
bits associated to sites in hexagons touching that interface,
only as needed to continue with the determination of the interface.
Note that the knowledge of $\beta$ suffices to determine if there
is an open crossing from the right boundary of $D$ above $p_0$ to
the left boundary of $D$: there is such a crossing if and only if
$\beta$ terminates on the left boundary of $D$.

Now let $\zeta'$ be the union of the bottom and left boundaries
of $D$, and let $\beta'$ be the interface from $p_0$ to $\zeta'$
that corresponds to the configuration $\omega'$ obtained by
flipping all the colors of the hexagons in $\closure {D}$
(alternatively, $\beta'$ is an interface that has black on the right
and white on the left).
Then $\beta'$ determines the existence of an open crossing from the
right boundary below $p_0$ to the left boundary.
Consequently, after the algorithm examines $\beta$ and $\beta'$,
the correct value of $Q$ is determined.
We now need to bound the revealment of this algorithm.

\begin{lemma}\label{l.Hdist}
Let $\p_W$ denote the counterclockwise arc of $\p D\setminus\zeta$
from $p_0$ to $\zeta$ (the white arc), and let
$\p_B$ denote the clockwise arc of
$\p D\setminus \zeta$ from $p_0$ to $\zeta$.
Let $H$ be a grid hexagon in $\closure D$.
Set $r_1:=\min\{\dist(H,\p_W),\dist(H,\p_B)\}$
and $r_2:=\max\{\dist(H,\p_W),\dist(H,\p_B)\}$.
Then
$$
\Pb{\p H\cap\beta\ne\emptyset \md p_0}
\le \alpha_2(r_1)\,\alpha^+\bigl(2\,r_1+O(1),r_2-r_1-O(1)\bigr)\,.
$$
\end{lemma}
\proof
Suppose that $\p H\cap\beta\ne\emptyset$.
Then there is a path in black hexagons from $\p H$
to $\p_B$ and there is a path in white
hexagons from $\p H$ to $\p_W$,
because the chains of hexagons along the two
sides of $\beta$ provide such paths.
Suppose, for example, that
$r_1=\dist(H,\p_W)$.
Let $w$ be a closest point to $H$ in $\p_W$.
Let $M$ be a half space that contains $D$,
which satisfies $\dist(\p M,w)\le C$ where $C=O(1)$.
(Here, we are using the fact that $D$ approximates a convex domain.)
Let $w'$ be a point closest to $H$ on $\p M$.
Then $\dist(w',H)\le r_1+C$.
Set $R_1=2\,r_1+2\,C+2\diam(H)$ and $ R_2=r_2-r_1-C-\diam(H)$,
and assume for now that $R_2>R_1$.
Consider next the annulus centered at
$w'$ with inner radius $R_1$ and outer radius $ R_2$.
Now, $M$ intersected with this annulus
contains a black crossing between the two boundary circles of this
annulus (because there is a black crossing from
$H$ to $\p_B$),
and there are black and white crossings between
$H$ and the circle of radius $r_1$ around the center
of $H$.  These events are independent given $p_0$,
which implies
the lemma in the case $r_1=\dist(H,\p_W)$, $R_2>R_1$.
If $R_2\le R_1$ and
$r_1=\dist(H,\p_W)$, we only need to consider the
crossings between $H$ and the circle of radius
$r_1$ around its center.
The case
$r_1=\dist(H,\p_B)$ is treated similarly.
\QED

\proofcont{Theorem~\ref{t.square}}
Fix a hexagon $H\subset\closure D$.
Note that the value of $r_1$ in the lemma does not depend on
$p_0$, since $r_1=\dist(H,\p D\setminus\zeta)$.
Let $z_0$ be the closest point to $H$ on the right boundary of
the square which $D$ approximates.
Observe that $r_2\ge |p_0-z_0|-O(1)$.
This implies that for every $r\in[1,R]$, $\Pb{r/2\le r_2<r}\le O(r/R)$.
Using the monotonicity of $\alpha^+$, Lemma~\ref{l.Hdist} therefore gives
\begin{multline*}
\Pb{\p H\cap\beta\ne\emptyset }
\\
\le O(1)\, \max_{0\le r_1\le R}\Bigl(
\alpha_2(r_1)
\sum_{j=0}^{\lceil \log_2 R\rceil}2^{-j}
\alpha^+\bigl(2\,r_1+O(1),2^{-j}R-r_1-O(1)\bigr)
\Bigr).
\end{multline*}
The same estimate also applies to $\beta'$, by symmetry.
Consequently,~\eref{e.expostrong} gives
\begin{multline*}
\Pb{H\text{ examined by algorithm}}
\\
\le R^{o(1)}\max_{0\le r_1\le R}\Bigl(
{r_1}^{-1/4}
\sum_{j=0}^{\lceil \log_2 R\rceil}2^{-j}
(2^{-j}R/r_1)^{-1/3}
\Bigr)
\\
\le R^{o(1)}\max_{0\le r_1\le R}\bigl( {r_1}^{1/12} R^{-1/3} \bigr)
= R^{-1/4+o(1)}\,,
\end{multline*}
as $R\to\infty$.
This proves the theorem in the case of the triangular grid.

The proof for the square grid is essentially the same.
Since in that case, we cannot use the values of the
critical exponents, we just use the
bounds $\alpha^+(r,r')\le C\,(r/r')^{\eps_0}$
and $\alpha_2(r_1)\le C\, r_1^{-\eps_0}$, which are
valid for some constants $C,\eps_0>0$.
The theorem follows.
\QED

We now give the proof of Theorem~\ref{th:crossingquant}. We will not prove
Theorem~\ref{th:crossingquantsquare}, but rather simply say that it is 
proved in a similar way.

\proofof{Theorem~\ref{th:crossingquant}}
Fix $\gamma < 1/8$.
Let $f$ be the indicator function $f(\omega)=1_{A_m'}(\omega)$.
By (\ref{e.page14}), we have
$$
N(A_m',m^{-\gamma})= \sum_{\emptyset\neq S \subseteq [n_m]} \hat{f}(S)^2 (1-2\,m^{-\gamma})^{2|S|}
$$
where $n_m$ is the number of sites in $D$.  By Corollary~\ref{c.square},
 with $\epsilon >0$ chosen so that $2\,\gamma+\eps < 1/4$,
this is at most
\begin{multline*}
C
\sum_{k=1}^{n_m} (1-2m^{-\gamma})^{2k}\, k\,m^{-1/4+\eps}\\
\le C\,{m^{-1/4+\epsilon}}
\sum_{k=1}^{\infty}k\, (1-2\,m^{-\gamma})^{\frac{m^\gamma}{2}\frac{4k}{m^\gamma}} 
\le C\,{m^{-1/4+\epsilon}}
\sum_{k=1}^{\infty} k\, e^{-4k{m^{-\gamma}}} 
\end{multline*}
It is easy to check that
$$
\sum_{k=1}^{\infty} k\, e^{-4k{m^{-\gamma}}} \le O(m^{2\gamma})\,,
$$
and so the result follows immediately.
\QED

\subsection{Annulus case}\label{ss.ann}

For the proof of Theorem~\ref{th:except}, we will need 
the following variant of Theorem~\ref{t.square}
regarding the percolation crossings of an annulus.

\begin{theorem}\label{t.annulus}
Let $2\le r< R$. 
Let $\croRr$ be the indicator function for the 
event that there is a crossing of the annulus
$\{z\in\R^2: r\le z\le R\}$
from the inner circle to the outer circle
by a cluster of white hexagons.
Then there is a randomized algorithm determining $\croRr$ with
\begin{equation}\label{e.td}
\delta\le r^{o(1)}\,\alpha(r,R)\,\alpha_2(r)\,.
\end{equation}
\end{theorem}

We stress that the $r^{o(1)}$ factor depends on $r$ only and not on $R$.
It is possible to replace the factor $r^{o(1)}$ by $O(1)$,
but in order to do this it seems that one would first need to
appeal to the analogue of~\eref{e.multip} for $\alpha_2$, which
is proved in the appendix.
In order to have a more direct proof of our main theorem,
we prefer, at this point, not to rely on the appendix.
By~\eref{e.onearmstrong} and~\eref{e.expo},
the right hand side in~\eref{e.td} is equal to
$R^{-5/48+o(1)}\,r^{-7/48}$, but its writing in~\eref{e.td} is more
suggestive and more useful.

\medskip

Since $\|\croRr\|^2=\alpha(r,R)$, 
Theorems~\ref{t.noise} and~\ref{t.annulus} give 

\begin{corollary}\label{c.croRr}
$$
\sum_{|S|=k}\hat\croRr(S)^2
\le 
k\,r^{o(1)}\,\alpha(r,R)^2\,\alpha_2(r)
$$
holds when $1\le r<R<\infty$ and $k>0$.
\QED
\end{corollary}

Before we prove Theorem~\ref{t.annulus},
we have to discuss the kind of interface that is used by the algorithm,
as it is slightly different from the interface used to determine
a possible crossing of a square.

\begin{figure}
\SetLabels
(.503*.50)$q_0$\\
\T(.44*.01)$p_0$\\
\B(.95*.61)$\beta^T$\\
\B(.95*.43)$\beta\setminus\beta^T$\\
\endSetLabels
\centerline{\epsfysize=2.5in%
\AffixLabels{%
\epsfbox{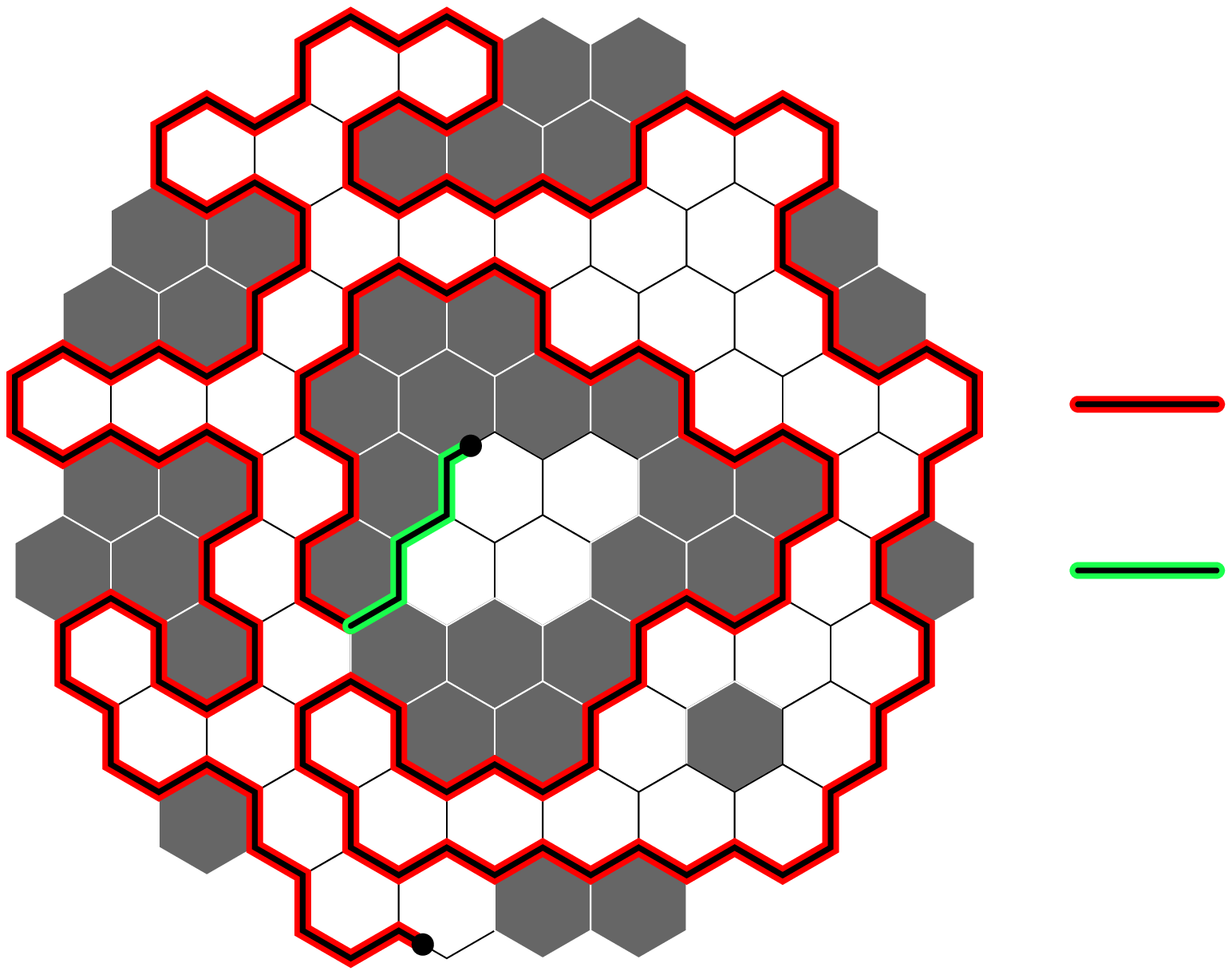}}%
}
\begin{caption} {\label{f.annulusinterface}The radial interface
$\beta$.}
\end{caption}
\end{figure}

Fix $R>0$ large.
Let $\closure D=\closure D_R$
be the union of the hexagons of the hexagonal grid 
that intersect the disk $|z|\le R$.
Let $V^*=V^*_R$ denote the set of vertices of the hexagonal grid
that are in $\closure D$.
Let $p_0$ be some point in $\p \closure D\setminus V^*$.
Let $H_0$ denote the hexagon containing the origin, and
let $q_0$ be some point in $\p H_0\setminus V^*$.
We define the {\bf radial interface} $\beta=\beta(R,p_0,q_0,\omega)$,
inductively as a simple path from $p_0$ to $q_0$.
(See Figure~\ref{f.annulusinterface}.)
The construction is segment by segment, and the concatenation
of the first $m$ segments will be denoted $\beta_m$.
If the (unique) hexagon in $\closure D$
containing $p_0$ is white [respectively, black], then 
the first segment $\beta_1$ of $\beta$ traverses the boundary
of that hexagon [counter-] clockwise 
until the first encounter with a point in $V^*$.
Suppose inductively $\beta_m$ has been constructed,
that it is a simple path, that
$p_m\in V^*$ and $p_0$ are the two
endpoints of $\beta_m$,
and that there is a path $\alpha$ in
the hexagonal grid in $\closure D$ from
$p_m$ to $q_0$ whose only intersection with $\beta_m$ is $p_m$.
The first step of such a path $\alpha$ is along an edge $\hat e$
starting at $p_m$. If there is just one possible $\hat e$ among all
such $\alpha$, then $\beta_{m+1}$ also uses that
edge $\hat e$. Clearly, there are at most two possible $\hat e$,
since the edge terminating at $p_m$ and used by $\beta_m$
cannot be used. If there are two possible $\hat e$, then
$\beta_{m+1}$ chooses between them according to the color of
the hexagon containing them both; i.e., the hexagon just encountered
by $\beta_m$. If that hexagon is white [respectively, black],
then the edge chosen is the one that traverses $H$
[counter-] clockwise. If the edge chosen contains $q_0$,
then the path stops at $q_0$ and the construction terminates.
Otherwise, $\beta_{m+1}$ is defined as the union of $\beta_m$ and
the chosen continuation edge. This completes the definition of $\beta$.

It is not hard to verify that for every simple path $\hat\beta$
in the hexagonal grid from $p_0$ to $q_0$ that stays in
$\closure D$, the probability that $\beta=\hat\beta$
is precisely $2^{-n}$ if $n$ is the number of
hexagons in $\closure D$ that intersect $\hat\beta$.
However, we will not use this fact.

Let $r\in[0,R]$.
We now define a truncated version of $\beta$, which will
suffice, as we will see, to determine $\croRr$.
We say that $\beta$ completed a [counter-] clockwise loop
at some dual vertex $v\in V^*$ if $v$ is visited by $\beta$
and there is a hexagon $H$ containing $v$ and  another  point
$u\in\p H$, which was visited by $\beta$ prior to
$v$, 
and the oriented arc of $\beta$ from $u$ to $v$ together with the
line segment $[v,u]\subset H$ form a [counter-] clockwise
loop surrounding $0$.
Let $\beta^T$ denote the initial segment of 
$\beta$ up to the first time in which $\beta$ completed
a counterclockwise loop around $0$ or until it hits
$q_0$, if there is no such loop.

\begin{lemma}\label{l.betadet}
The truncated radial interface $\beta^T$ meets 
the disk $|z|\le r$ if and only if $\croRr=1$.
\end{lemma}
\proof
Let $[u,v]$ be an edge in $\beta^T$, with $u$ occuring
before $v$ along $\beta$.
We claim that if the hexagon $H$ immediately to the right
of $[u,v]$ is contained in $\closure D$, then it is white.
Indeed, suppose the contrary. Let $w$ be the first vertex along 
$\beta$ in which $\p H$ is visited, and let
$\beta^w$ be the initial segment of $\beta$ from $p_0$
to $w$.  Observe that the counterclockwise arc from 
$w$ to $v$ is a feasible continuation of $\beta^w$,
since $\beta$ contains a path from $v$ to $q_0$ and
there is no other point but $w$ in
$\p H\cap\beta^w$.
Since we are assuming that $H\subset\closure D$ is black,
it follows that the immediate continuation of $\beta^w$
was along $\p H$ in the counterclockwise direction until
some $w'\in V^*$ is hit.

Consider the directed cycle obtained by joining the line segment
$[u,w']$ to the arc of $\beta$ from $w'$ to $u$.
This directed cycle surrounds $v$, because $w$ is connected
in $\beta^w$ to $p_0$, which is certainly in the unbounded
component of this cycle, and the line segment
$[v,w]$ intersects the cycle precisely once,
crossing the line segment $[u,w']$ inside $H$.
Moreover, if we consider the orientation in which
these two line segments cross, we conclude that
the cycle surrounds $v$ counterclockwise.
Because the arc of $\beta$ from $v$ to $q_0$
does not intersect the cycle, we conclude that the cycle also
surrounds $0$ counterclockwise.
This contradicts the definition of the truncated path $\beta^T$,
since we are assuming $v\in\beta^T$.
This verifies our claim, that to the right of edges in
$\beta^T$ are only white hexagons and hexagons that are not
contained in $\closure D$.

Note that if $e$ and $e'$ are two consecutive segments
along $\beta$, then the hexagon to the right of $e$ is
either the same as the one to the right of $e'$, or these
hexagons are adjacent. We therefore conclude that every
white hexagon visited by $\beta^T$ is connected by
a chain of white hexagons to $\p D$.
Therefore, if $\beta^T$ hits the set $|z|\le r$,
then clearly $\croRr=1$.

Now suppose that $\beta^T$ does not hit $|z|\le r$.
This implies that $\beta^T$ has terminated by completing a
counterclockwise loop around the set $|z|\le r$.
Consider a hexagon $H$ on the inner boundary of this loop.
Because the orientation of the loop is
counterclockwise, the first time in which $H$ is
visited, the path chose to traverse $\p H$ counterclockwise.
This implies that the hexagon is black.
Thus, there is a loop of black hexagons in $\closure D$
that surrounds the set $|z|\le r$. This implies that $\croRr=0$.
\QED

We can now specify the algorithm promised by Theorem~\ref{t.annulus}.
The algorithm starts by selecting the point $p_0$ uniformly along $\p\closure D$,
and selecting $q_0$ arbitrarily in $\p H_0\setminus V^*$.
It then proceeds to inspect the colors of the hexagons necessary to develop
the truncated interface $\beta^T$, until it terminates or
hits the set $|z|\le r$. At that point, the correct value
of $\croRr$ is determined, by Lemma~\ref{l.betadet}.

In order to bound the revealment of this algorithm,
it will be convenient to introduce a different interface, which in the end will
turn out to be equivalent to $\beta$.

Let $\hat D$ denote the branched double cover of $\closure D$ about $0$,
and let $\phi:\hat D\to \closure D$ denote the projection map.
Concretely, define $\hat D$ as the preimage of
$\closure D$ under the map $\phi(z)=z^2$.
Let $\hat p_0$ be one of the preimages of $p_0$ under $\phi$,
and let $\hat q_0$ be one of the preimages of $q_0$.
Let $\hat H_0$ be the closure of one of the connected components
of $\phi^{-1}(H_0)\setminus [\hat q_0,-\hat q_0]$.
Let $\mfh$ denote the set of hexagons $H$ that are contained
in $\closure D$.
Let $\hmfh$ denote the set of connected components
of preimages $\phi^{-1}(H)$, $H\in\mfh$, except that
the single preimage of $H_0$ is replaced by
the two sets $\hat H_0$ and $-\hat H_0$.
Note that if $\hat H\in\hmfh$,
then $-\hat H\in\hmfh$ and $\phi(\hat H)=\phi(-\hat H)\in\mfh$.
Let $\hmfh' \subset \hmfh $ be a maximal collection
of elements of $\hmfh$ with the property that
$\hmfh'\cap \{-\hat H:\hat H\in\hmfh'\}=\emptyset$.
Now color at random each of the elements of $\hmfh'$ white or
black independently, with probability $1/2$.
For every $\hat H\in\hmfh'$, let $-\hat H$ have the opposite
color to the color of $\hat H$.

Now let $\hat \beta$ denote the chordal interface in $\hat D$
from $\hat p_0$ to $-\hat p_0$, with white cells on the right
and black cells on the left, as defined in the simply connected setting
in Subsection~\ref{ss.simply}.
That is, we consider the exterior of the counterclockwise arc from $\hat p_0$
to $-\hat p_0$ along $\p \hat D$ as white, the exterior of the complementary
arc as black, and take $\hat \beta$ as the interface between white and
black passing through $\hat p_0$ and through $-\hat p_0$.
Finally, let $\beta^\dagger:= \phi(\hat\beta)\setminus \interior(H_0)$.

\begin{lemma}
Given $p_0$ and $q_0$,
the laws of $\beta^\dagger$ and of $\beta$ are the same.
\end{lemma}
(In this statement, we consider $\beta$ as a set, and forget
about the fact that it has the structure of an oriented path.)

\proof
The map $z\mapsto -z$ preserves $\hat\beta$, by the symmetry
of the interface. Consequently, near every point $p\in \beta^\dagger\setminus \{p_0,q_0\}$,
$\beta^\dagger$ looks like a piecewise linear path. Moreover,
$\beta^\dagger$ is connected and contains $p_0$. Since a compact simple path has two endpoints,
we conclude that $q_0\in\beta^\dagger$ as well.
We now consider building $\beta^\dagger$ by adding one segment at a time.
When it hits a previously unvisited hexagon $H$ which is contained in $\closure D$,
it is equally likely (given its past) to turn right or left.
(This is because both preimages of $H$ are unvisited by both preimages of the
past of $\beta^\dagger$.) When it hits a previously visited hexagon (or
a hexagon that is not contained in $\closure D$), it turns in such a way
that it will eventually be able to reach $q_0$ without crossing itself,
and this uniquely specifies this turn. Consequently, the lemma follows.
\QED

\begin{remark}
The radial interface converges to radial SLE(6) as the mesh tends to zero.
\end{remark}

\proofof{Theorem~\ref{t.annulus}}
Given all of our preparations, the proof is rather easy.
We have shown that the above algorithm provides the correct answer.
It therefore remains to estimate its revealment.
Consider some hexagon $H\subset\closure D$.
We want to prove that the right hand side of~\eref{e.td}
is an upper bound for the probability that $H$ is examined.
Let $a:=\dist(0,H)$, $b:=\dist(H,\p\closure D)$ and $c:=\dist(p_0,H)$.
Let $S_1$ be the disk of radius $(a\wedge b)/2$ concentric with $H$,
and let $\hat S_1$ be one of the two connected components
of $\phi^{-1}(S_1)$.
We have to bound the probability that the algorithm inspects $H$.
Clearly, we may assume $a\ge r-O(1)$.
For $H$ to be inspected,
 $\beta^T$ has to get to the circle
$|z|=R\wedge (2\,a)$.  This probability
is $\alpha(2\,a,R)$.
Given that this has happened, how can we estimate the
probability that $\beta$ is adjacent to $H$?
At this point, we use the equivalence of
$\beta$ and $\beta^\dagger$.
The information that $\beta^\dagger$ reached the circle
$|z|=R\wedge (2\,a)$ bears no impact on the distribution
of the colors of the cells in $\hmfh$ whose images under
$\phi$ intersect $S_1$. (Here, we assume that $a$ is not too small,
so that the corresponding sets of cells are disjoint.
Certainly $a>10$ suffices. If $a$ is smaller,
then the estimate we are now striving for is trivial.)
Since there is no hexagon intersecting both $\hat S_1$ and $-\hat S_1$,
it follows that the conditional distribution of the
colors of the cells meeting $\hat S_1$ is uniform
i.i.d. Consequently, the conditional probability that
$\beta^\dagger$ hits $H$ is bounded by $\alpha_2((a\wedge b)/2)$.
Thus,
$$
\Pb{H\text{ visited}}\le
O(1)\,\alpha(2\,a,R)\,
\alpha_2((a\wedge b)/2)\,.
$$
In the case $b\ge a\ge 2\,r$, we may use independence on disjoint
sets to conclude that 
\begin{equation*}
\begin{aligned}
&\Pb{H\text{ visited}}
 \\
&\qquad
\le
O(1)\, \alpha_2(r)\,\alpha_2(2\,r,a/2)\,\alpha(2\,a,R)
\le
O(1)\,
\alpha_2(r)\,\alpha(2\,r,a/2)\,\alpha(2\,a,R)
\\
&\qquad
\overset{\eref{e.rsw}}{\le}
O(1)\,
\alpha_2(r)\,\alpha(r,2\,r)\,\alpha(2\,r,a/2)\,
\alpha(a/2,2\,a)\,\alpha(2\,a,R)
\\
&\qquad
\overset{\eref{e.multip}}{\le}
O(1)\,
\alpha_2(r)\,\alpha(r,R)\,.
\end{aligned}
\end{equation*}
On the other hand, if $b\ge a$ and $a\le 2\,r$,
then we use our assumption $a\ge r-O(1)$ and~\eref{e.expo}
to get
$\alpha_2(a/2)\le r^{-1/4+o(1)}\le\alpha_2(r)\,r^{o(1)}$,
which is also sufficient.

In the case $b< a$, a similar argument 
(and similar to the proof of Lemma~\ref{l.Hdist})
shows that
$$
\Pb{H\text{ visited }\md c}\le O(1)\,\alpha_2(b/2)\,\alpha^+\bl(2\,b+O(1),c-b-O(1)\br).
$$
Next, picking a constant $q\in(1/4,1/3)$, we then have by the above and~\eref{e.expostrong}
$$
\Pb{H\text{ visited }\md c}
\le O(1)\, \alpha_2(b/2)\, (c/b)^{-q}.
$$
As in the proof of Theorem~\ref{t.square}, we have $\Pb{2^j\le c<2^{j+1}} \le O(1)\,2^j/R$.
It easily follows that
$$
\Pb{H\text{ visited}}\le
O(1) \,\alpha_2(b/2)\, (R/b)^{-q}.
$$
If $b/2>r$, then we may estimate
$$
\alpha_2(b/2) \le \alpha_2(r)\, \alpha_2(2\,r,b/2) \le \alpha_2(r)\, \alpha(2\,r,b/2)
$$
and 
$$
(R/b)^{-q}
\overset{\eref{e.expostrong}}{\le}
O(1)\, \alpha(b/2,R)
$$
and we get from~\eref{e.multip} and the above
$$
\Pb{H\text{ visited}}\le O(1)\, \alpha_2(r)\, \alpha(r,R)\,.
$$
If $b/2\le r$, we use instead
$$ 
\alpha_2(b/2)\, (R/b)^{-q}\le O(1)\, \alpha_2(b/2)\, \alpha_2(b/2,r)\, \alpha(r,R)
\overset{\eref{e.expostrong}}{\le}
r^{o(1)}\, \alpha_2(r)\, \alpha(r,R)\,.
$$
This completes the proof.
\QED

\begin{remark}
It is easy to see that if we assume the analogue of~\eref{e.multip} for 
$\alpha_2$ proved in the appendix,
then the $r^{o(1)}$ term in~\eref{e.td} can be replaced by $O(1)$.
\end{remark}

\section{Exceptional times}\label{sec:except}

In this section we prove Theorem \ref{th:except}. 
We point out that the absolute key necessary step is to get a good
bound on the correlation for an event occurring at two different but
close by times. Once this is done, the rest is fairly standard.
Proposition A16 in Lawler~\cite{Lawler}
indicates this general type of argument.

\proofof{Theorem~\ref{th:except}}
By Kolmogorov's 0-1 law, it suffices to prove that with positive probability
there are times in $[0,1]$ when the origin is in an infinite cluster.
Fix $R>2$ large and let $V_{t,R}$ be the event
that at time $t$ there is an open path from the origin to distance $R$ away.
We then let
$$
X=X_R:= \int_0^1 1_{V_{t,R}} \, dt 
$$
be the Lebesgue measure of the set of times in $[0,1]$ at which $V_{t,R}$ occurs.
The first moment of $X$ is given by
$$
\Eb{X}=\int_0^1 \Pb{V_{t,R}}\,dt=\Pb{V_{0,R}} =\alpha(R)\,.
$$
The second moment is
\begin{equation}\label{e.2nd}
\Eb{X^2}=
\EB{\int_0^1\int_0^1 1_{V_{s,R}}\,1_{V_{s',R}}\,ds\,ds'}
=
\int_0^1\int_0^1 \Pb{V_{s,R}\cap V_{s',R}}\,ds\,ds'.
\end{equation}

For each site $v$ we let 
$$
\chi_v^s:=\begin{cases} -1&v \text{ is open at time } s\\
1&\text{otherwise},
\end{cases}
$$
and for a finite set of sites $S$ set
$$
\chi_S^s:=\prod_{v\in S}\chi_v^s\,.
$$
Fix $s,s'\in[0,1]$, and set $t:=|s-s'|$.
Recall that the state of a site $v$ flips between closed and open
with rate $1/2$. Equivalently, we may think of the state as being
re-randomized with rate $1$.
Consequently, $\Pb{\chi_v^{s'}=\chi_v^s \md \omega_s} = e^{-t}+(1-e^{-t})/2
=(1+e^{-t})/2$,
and hence,
$$
\Eb{\chi_v^s\,\chi_v^{s'}}=\exp(-t)
\,,\qquad
\Eb{\chi_S^{s}\,\chi_S^{s'}}=\prod_{v\in S}\Eb{\chi_v^s\,\chi_v^{s'}}=\exp(-t\,|S|)\,.
$$
Moreover, if $S\ne S'$, then $\Eb{\chi_S^s\,\chi_{S'}^{s'}}=0$.
Consequently, if $f$ is a function depending on the states of finitely many lattice
points and has the expansion
$f(\omega)=\sum_S\hat f(S)\,\chi_S(\omega)$, then
\begin{equation}\label{e.timeexpansion}
\Eb{f(\omega_s)\,f(\omega_{s'})}=
\sum_S\hat f(S)^2 \exp(-t\,|S|)\,.
\end{equation}

Let $f_r^R(\omega)$ be as in Theorem~\ref{t.annulus}.
Fix some $t\in (0,1]$ and let $r\in [2,R)$.
Clearly, $0\le f_0^R(\omega)\le f_0^r(\omega)\,f_{2r}^R(\omega)$ for every $\omega$.
Consequently,
\begin{multline*}
\Pb{\Vss}=
\Eb{f_0^R(\omega_s)\,f_0^R(\omega_{s'})}
\le
\Eb{f_0^r(\omega_s)\,f_{2r}^R(\omega_s)\,
f_0^r(\omega_{s'})\,f_{2r}^R(\omega_{s'})}
\\
=
\Eb{f_0^r(\omega_s) \, f_0^r(\omega_{s'})}
\Eb{ f_{2r}^R(\omega_s)\,f_{2r}^R(\omega_{s'})}
\le
\Eb{f_0^r(\omega_s) } \Eb{ f_{2r}^R(\omega_s)\,f_{2r}^R(\omega_{s'})} .
\end{multline*}
(To obtain the second equality, we have used the independence on 
disjoint sets of sites.) Thus,
$$
\Pb{\Vss}\le
\alpha(r)\,
\Eb{ f_{2r}^R(\omega_s)\,f_{2r}^R(\omega_{s'})}
=\alpha(r)\, \sum_S e^{-t|S|}\hat f_{2r}^R(S)^2.
$$
The latter sum restricted to $S$ with $|S|=k$ for fixed $k\ne 0$
is estimated using Corollary~\ref{c.croRr},
while for $k=0$, we use $\hat f_{2r}^R(\emptyset)=\alpha(2r,R)$. This yields
$$
\Pb{\Vss}\le
\alpha(r)\,\Bigl(\alpha(2\,r,R)^2+
r^{o(1)}\,\sum_{k=1}^\infty 
e^{-kt}
k\,\alpha(2\,r,R)^2\, \alpha_2(r)\Bigr).
$$
It is easy to check that $\sum_{k=1}^\infty k\,e^{-kt}\le O(t^{-2})$.
This and the inequalities~\eref{e.multip} and~\eref{e.rsw} allow
us to write this estimate as
\begin{equation}\label{e.zz}
\Pb{\Vss}\le
O(1)\,
\alpha(R)^2\,\alpha(r)^{-1}\bl(1+r^{o(1)}\,t^{-2}\,\alpha_2(r)\br)\,.
\end{equation}
We proved the above claim for all 
$t\in (0,1]$ and $r\in [0,R)$ but
now we observe that~\eref{e.zz} is also trivially true
when $r \ge R$ as well. 
We now choose $r=2\,t^{-8}=2\,|s-s'|^{-8}$.
Applying this in~\eref{e.zz} with~\eref{e.onearm} and~\eref{e.expo} gives
\begin{equation}\label{e.uselater}
\Pb{V_{s,R}\cap V_{s',R}}
\le O(1)\,\alpha(R)^2\, |s-s'|^{-5/6+o(1)}.
\end{equation}
Hence,
\begin{equation}\label{e.fin}
\int_0^1\int_0^1\Pb{V_{s,R}\cap V_{s',R}}\,ds\,ds'
\le O(1)\,\alpha(R)^2.
\end{equation}

The Cauchy-Schwarz inequality tells us that
$$
\Pb{X>0} \ge 
\frac{
\Eb{X}^2
}{
\Eb{X^2}
}\,.
$$
Consequently, the above inequality, the fact that $\Eb{X}=\alpha(R)$,
the expression~\eref{e.2nd} for $\Eb{X^2}$ and~\eref{e.fin}
show that $\inf_{R>0}\Pb{X_R>0}>0$.
Let $T_R:=\{t\in[0,1]:V_{t,R}\text{ holds}\}$.
 Fatou's lemma tells us that with positive probability
$T_R\ne \emptyset$ for infinitely many $R\in\N$.
Since $T_R\supset T_{R'}$ when $R'>R$,
this implies that 
$$
\Pb{\cap_{R>0}\,\, \{T_R\ne\emptyset\}}>0.
$$
Our goal is to show that $\Pb{\bigcap_R T_R \neq \emptyset}>0$.
Since the $T_R$'s are not closed sets,
$\cap_{R>0}\,\, \{T_R\ne\emptyset\}$
does not immediately imply $\bigcap_R T_R\ne\emptyset$.
(The reason that $T_R$ is not necessarily closed is that the set of times at
which an edge is open is not a closed set since we have a right continuous
process.) This technicality is taken care of by the following lemma from \cite{HPS}. 

\begin{lemma}\label{l.HPS} (\cite{HPS})
Let $0<p<1$ and let $G$ be any graph where $\pi_p(\calC)=0$.
Let $\{\omega_t\}$ represent our dynamical percolation process in
that $\omega_t(v)$ is the state of vertex $v$ at time $t$.
Consider the process $\{\bar\omega_t\}$ obtained from $\{\omega_t\}$
by setting, for every vertex $v$,
the set  $\{t\, : \, \bar\omega_t(v) =1\}$ to be 
the closure of the set $\{t\, : \, \omega_t(v) =1\}$.
Then $\bPsi_p$-a.s., for every vertex $v$ we have
$$
  \{t \in [0, \infty) \: :\:  v \, 
\mbox{ \rm percolates in } \bar\omega_t \} \, = \,
 \{t \in [0, \infty) \: :\:  v \, \mbox{ \rm percolates in } \omega_t \} \,.
$$
In particular, a.s.\ this set of times is closed.
\end{lemma}

Returning to our proof,
let $\overline{T_R}$ be the closure of $T_R$. It is easily checked that
$$
\bigcap_{R>0}\overline{T_R}= 
  \{t \in [0,1] \: :\:  0 \, 
\mbox{ \rm percolates in } \bar\omega_t \},
$$
where $\{\bar\omega_t \}$ is defined in Lemma \ref{l.HPS}.
By compactness, if the $T_R$'s are all nonempty, it follows that 
$\bigcap_R\overline{T_R}$ is nonempty. This implies that there 
is some time at which $\bar\omega_t$ percolates and hence by
Lemma \ref{l.HPS}, some time at which the original process $\omega_t$ percolates.
\QED

For future reference, we note that Lemma~\ref{l.HPS} implies that a.s.
\begin{equation}\label{e.ctr}
{\bigcap_{R>0}\overline{T_R}= 
\bigcap_{R>0}{T_R} }.
\end{equation}

\section{Hausdorff dimension of exceptional times} \label{sec:hd}

In this section, we prove Theorem \ref{th.hd}. This is separated into two theorems,
Theorem \ref{t.lowerbound} and Theorem \ref{t.upperbound}, where lower and upper bounds
are given. We point out however that the lower bound is simply a refinement of
the argument for proving that there exist exceptional times.
First note that the fact that the  Hausdorff dimension is an almost sure
constant follows immediately from ergodicity.

\begin{theorem}\label{t.lowerbound}
A.s.,  the Hausdorff dimension of the set of exceptional times is at least $1/6$.
\end{theorem}

\proof
Fix $\gamma < 1/6$. It suffices by ergodicity and countable additivity
to show that with positive probability,
the set of exceptional times in $[0,1]$ at which the origin percolates
has Hausdorff dimension at least $\gamma$.
For each integer $R$, let as before $V_{t,R}$ be the event
that at time $t$ there is a path from the origin to distance $R$ away and
define a random measure $\sigma_R$ on $[0,1]$ by
$$
\sigma_R(S)=\frac{1}{\alpha(R)}\int_{S}1_{V_{t,R}}dt
$$
for each Borel set $S\subset [0,1]$.

The results in the previous section immediately give that
$\E[\| \sigma_R \|] =1$ and $\E[\| \sigma_R \|^2] \le O(1)$
where $\| \sigma_R \|$ denotes the total variation of the measure $\sigma_R$.

Cauchy-Schwarz gives
$$
\Eb{\|\sigma_R\|^2}^{1/2}\,\Pb{\|\sigma_R\|>1/2}^{1/2}
\ge
\Eb{\|\sigma_R\| 1_{\|\sigma_R\|>1/2}}
\ge
\Es{\|\sigma_R\| }-1/2=1/2.
$$
Consequently, 
$\Pb{\| \sigma_R \| > 1/2} \ge  C_1 $
for some constant $C_1 >0$.
Given a measure $m$ on $[0,1]$ and $\gamma >0$, let
$$
{\cal E}_{\gamma} (m) = \int \!\!\int {|t-s|}^{- \gamma} \, dm (t) \,   dm(s) .
$$
Note that
$$
\Eb{{\cal E}_{\gamma} (\sigma_R)}
=
\EB
{\int_0^1\!\int_0^1\frac{d\sigma_R (t) \,   d\sigma_R(s)}{ {|t-s|}^ \gamma}  }=
\int_0^1\!\int_0^1
\frac{\Pb{{V_{t,R}}\cap{V_{s,R}} }}{{\alpha(R)^2}\, |t-s|^\gamma}\, dt\,ds \,.
$$
Therefore, 
by~\eref{e.uselater} and $\gamma<1/6$,
$$
C_2:=\sup_R
\Eb{{\cal E}_{\gamma} (\sigma_R)}<\infty\,.
$$

By Markov's inequality, for all $R$ and for all $T$,
$$
\Pb{{\cal E}_{\gamma} (\sigma_R)\ge  C_2T} \le 1/T.
$$
Choose $T$ so that $1/T < C_1/2$.
Letting 
$$
U_R=\{\| \sigma_R \| > 1/2\}\cap \{{\cal E}_{\gamma} (\sigma_R)\le  C_2T\},
$$
by the choice of $T$, we have that 
$$
\Pb{U_R} \ge C_1/2.
$$
By Fatou's lemma,
$$
\Pb{\limsup_{R\to\infty} U_R} \ge C_1/2.
$$
We now show that on the event $\limsup_R U_R$, the Hausdorff dimension of the 
set of percolating times in $[0,1]$ is at least $\gamma$. 
Let $\overline{T_R}$ again be the closure of the set of times in $[0,1]$ at which 
there is a path from the origin to distance $R$ away. Clearly $\sigma_R$ is supported on
$\overline{T_R}$.
By~\eref{e.ctr}, it suffices to prove that
$\bigcap_{R>0}\overline{T_R}$ has Hausdorff dimension at least $\gamma$
on the event $\limsup_R U_R$.
This is achieved in the following
(deterministic) lemma, which completes the proof.
\QED

\begin {lemma}\label{l.deterministic}
Let $D_1\supseteq D_2 \supseteq D_3 \ldots$
be a decreasing
sequence of compact subsets of $[0,1]$,
and let $\mu_1,\mu_2,\dots$ be a sequence of positive measures
with $\mu_n$ supported on $D_n$.
Suppose that there is a constant $C$ such that for infinitely many
values of $n$, we have
\begin{equation}\label{e.2cond}
\| \mu_n \| > 1/C,\quad \text{ and } \quad{\cal E}_{\gamma} (\mu_n)\le C.
\end{equation}
Then the Hausdorff dimension of $\bigcap_n D_n$ is at least $\gamma$. 
\end {lemma}

\proof
Choose a sequence of integers $\{n_k\}$ for which
(\ref{e.2cond}) holds.
Note that $\|\mu_{n_k}\|^2\le \mathcal E_\gamma(\mu_{n_k})\le C$.
 By compactness, 
choose a further subsequence $\{n'_{k}\}$ of $\{n_k\}$ so that
$\mu_{n'_{k}}$ converges weakly to some positive measure $\mu_{\infty}$.
Clearly $\mu_{\infty}$ is supported on $\bigcap_n D_n$ and
$\| \mu_\infty \| \ge 1/C$. For all $M$, we have that
\begin{multline*}
\int \!\!\int {|x-y|}^{- \gamma} \wedge M \, d\mu_\infty (x) \,   d\mu_\infty(y) =
\\ \lim_{k\to\infty}
\int \!\!\int {|x-y|}^{- \gamma} \wedge M \, d\mu_{n'_k} (x) \,   d\mu_{n'_k}(y) \le C.
\end{multline*}
Now let $M\to\infty$ and apply the monotone convergence theorem to conclude that
$$
\int \!\!\int {|x-y|}^{- \gamma}\, d\mu_\infty (x) \,   d\mu_\infty(y) 
\le C.
$$
Since $\| \mu_\infty \| > 0$, it now follows from Frostman's theorem
(see for example,~\cite{Kahane})
that the Hausdorff dimension of $\bigcap_n D_n$ is at least $\gamma$. 
\QED

\begin{theorem}\label{t.upperbound}
A.s., 
the Hausdorff dimension of the set of exceptional times is at most $31/36$.
\end{theorem}

\proof
Let $U_n$ be the event that there is a time in $[0,1/n]$ for which
the origin percolates. Since the set of vertices which are open for some $t\in [0,1/n]$
is an i.i.d.\ process with density $1/2 +(1-e^{-1/(2n)})/2\le 1/2+1/n$, it is immedate that
$$
\Pb{U_n}\le \pi_{\frac{1}{2}+\frac{1}{n}} (\calC_0) ,
$$
where $\calC_0$ is the event that the origin percolates.
By page 3 of \cite{SW}, for every $\epsilon >0$, there is a $C$ so that
\begin{equation}\label{e.offcrit}
\pi_{\frac{1}{2}+\frac{1}{n}} (\calC_0)\le 
C\, n^{\eps-5/36}.
\end{equation}
Now let 
$$
N_n=\sum_{j=1}^n 1_{U_{j,n}}\,,
$$
where $U_{j,n}$ is the event that there is a time in $[(j-1)/n,j/n]$ for which
the origin percolates (so $U_{1,n}=U_n$ above). 
By the above, we have that $\Eb{N_n}\le C n^{\frac{31}{36}+\epsilon}$.
It follows that
$$
\lim_n \frac{\Eb{N_n}}{n^{\frac{31}{36}+2\epsilon}}=0
$$
and so from Fatou's lemma, we get
$$
\EB{\liminf_n \frac{N_n}{n^{\frac{31}{36}+2\epsilon}}}=0.
$$
Therefore
$$
\liminf_n \frac{N_n}{n^{\frac{31}{36}+2\epsilon}}=0
$$
a.s. This says that a.s.\ for infinitely many $n$, the set of exceptional
times in $[0,1]$ at which the origin percolates
can be covered by $n^{\frac{31}{36}+2\epsilon}$ intervals of length $1/n$.
Hence, the Hausdorff dimension of the set of these exceptional times is 
at most $\frac{31}{36}+2\epsilon$ a.s. By countable additivity, we are done.
\QED

\begin{remark}
The upper bound will be proved again by a different argument when
we prove Theorem~\eref{th:cones}.
The above proof is included 
here, because it is shorter.
 One should nonetheless
point out that the above argument uses~\eref{e.offcrit},
while the argument below is more self contained.
\end{remark}

\section{Exceptional times for $k$-arm events} \label{sec:wedgeslower}

In this section, we give the proofs of the lower bounds in
Theorems \ref{th:wedges} and \ref{th:cones},
but generally omit those details which are the same as in the
corresponding proofs of Theorems~\ref{th:except} and~\ref{t.lowerbound}.

For $\theta>0$ and integer $k\ge 1$, let
$A_{W_\theta}^k(r,R)$ be the event that
we have $k$ disjoint crossings of alternating colors (with black most clockwise)
between distances $r$ and $R$ of the origin in $W_\theta$
and let $\alpha_{W_\theta}^k(r,R)=\Pb{A_{W_\theta}^k(r,R)}$. 
If $r$ is suppressed, then it is taken to be $10\,k$.

We will, of course, need the asymptotics of $\alpha_{W_\theta}^k(r,R)$.
For this purpose, conformal invariance will be used.
Although when $\theta>2\pi$ the surface
$W_\theta$ is not planar and conformal invariance is usually stated for planar domains,
the proof of conformal invariance certainly holds in this setting.
The asymptotic decay as $R/r\to\infty$ of the probability of $k$ disjoint crossings
between distances $r$ and $R$ in $W_{\theta}$ from the origin in 
the percolation scaling limit is determined using conformal invariance.
Specifically, the map $z\mapsto z^{\pi/\theta}$ maps $W_{\theta}$ to the
upper half plane, and we may conclude from~\eref{e.halfpk}
that the decay (for the percolation
scaling limit)  is of the form $(R/r)^{\frac{-\pi k(k+1)}{6\theta}+o(1)}$,
as $R/r\to\infty$ while $k$ stays fixed.
Then, one can conclude, as for the other exponents we have discussed,
that for $R\ge r\ge 10\,k$,
\begin{equation}\label{e.awt}
\alpha_{W_\theta}^k(r,R)= (R/r)^{\frac{-\pi k(k+1)}{6\,\theta}+o(1)},
\end{equation}
as $R/r\to\infty$ while $k$ is fixed, using the argument given in~\cite{SW}. 
We will also use the fact that the quasi-multiplicativity relations~\eref{e.multip}
and~\eref{e.rsw} hold for $\alpha_{W_\theta}^k$ and
for $\alpha_2$. This is proved in the appendix;
see Remark~\ref{r.qmc}.

\proofof{lower bound in Theorem~\ref{th:wedges}}
We first handle the case $k=1$ and therefore abbreviate temporarily
$\alpha_{W_\theta}^1(r,R)$ by $\alpha_{W_\theta}(r,R)$.
(A different approach will be needed for $k\ge 2$.)
Let $D$ be the union of the hexagons in $W_\theta$ that contain
points whose distance from the origin is in $[r,R]$.
Let $\p_R D$ and $\p _r D$ denote the set of points in
$\p D$ that are at distance $\ge R$ [respectively,
$\le r$] from the origin.
Also, we denote by $\p^0 D$ and $\p^\theta D$,
the components of $\p D\cap\p W_\theta$ that are at angle
about $0$ [respectively, about $\theta$] in radial coordinates on $W_\theta$.

The algorithm we use to determine if
there exists a crossing of $D$ is essentially the same as
the algorithm determining the existence of a left to right crossing of a square,
where $\p_R D$
plays the role of the right side of the square and
where $\p_r D$
plays the role of the left side of the square.
(This is of course crucial; if we reversed things, then the hexagons
near the inner circle would be revealed with too high of a probability.)
It is clear that this algorithm works and so we now need
to compute its revealment.
We will show that the revealment is
\begin{equation}\label{e.revealment}
O(1)\,\alpha_2(r)\,\alpha_{W_\theta}(r,R)\,.
\end{equation}
Using~\eref{e.awt} and~\eref{e.expostrong}, one can show that this is
essentially (i.e., up to some $O(1)$ factor) monotone decreasing
in $r$ in the relevant range $\theta>8\,\pi/3$.

We just look at the first interface arising in the algorithm,
the one which terminates when it hits $\p_r D\cup\p^\theta D$,
since the estimates for the second interface will be essentially the same.

Fix some hexagon $H\subset D$.
Let $s=\dist(H,\p_R D\cup \p^0 D\cup \{0\})$ with $0$ denoting the origin.
We also use $|p|$ to denote distance from $0$ in $W_\theta$.
We distinguish several different cases.

\medskip\noindent
{\bf Case 1}: $\dist(H,\{0\})=s$. 
For $H$ to be visited, we need our 2 arms event holding within distance $s/2$ of $H$
and a crossing of the desired color between distance $2s$ and distance $R$ from the 
origin. These are independent and we get that $H$ is visited with probability 
at most $\alpha_2(s/2)\,\alpha_{W_\theta}(2s,R)$.
By the analogues of~\eref{e.multip} and~\eref{e.rsw} for
$\alpha_2$ and $\alpha_{W_\theta}$, this is compatible with
our claimed revealment~\eref{e.revealment}.

\medskip\noindent
{\bf Case 2}: $\dist(H,\p_R D)=s$. 
As in the proof of Theorem~\ref{t.annulus} with $c:=\dist(p_0,H)\wedge (R/2)$,
we obtain
$$
\Pb{H\text{ visited }\md p_0}\le O(1)\,\alpha^+(2s+O(1),c-s-O(1))\,\alpha_2(s/2).
$$
Proceeding as in that proof, this is also compatible with
our claimed revealment~\eref{e.revealment}.

\medskip\noindent
{\bf Case 3}: $\dist(H,\p^0 D)=s$. 
Let $w\in \p^0 D$ so that $\dist(H,w)=s$. We separate Case 3 into 3 subcases.

\medskip\noindent
{Case 3(a)}: $s \ge |w|/2$.
Then the triangle inequality gives $\dist(H,0)\le 3s$. 
For $H$ to be visited, we need our 2 arms event holding within distance $s/2$ of 
$H$ and a path of the desired type between distance $4s$ and distance $R$ 
from the origin. These are independent and we get that $H$ is visited with 
probability at most $\alpha_2(s/2)\,\alpha_{W_\theta}(3s,R)$,
which is compatible with our claimed revealment~\ref{e.revealment}.

\medskip\noindent
{Case 3(b)}: $s \le |w|/2 \le R/4$.
For $H$ to be visited, we need our 2 arms event holding within distance $s/2$ of $H$,
a white crossing in the half annulus centered at $w$ with outer radius $|w|$ and inner radius $2s$
(which is identical to a half-annulus in a half-plane;
$|w|\le R/2$ guarantees that the above half-annulus does not intersect $\p_R D$)
and a white crossing between distance $2|w|$ and distance $R$ from the 
origin.
These are independent and we
get that $H$ is visited with probability at most 
$$
\alpha_2(s/2)\,\alpha_{+}(2s,|w|)\, \alpha_{W_\theta}(2|w|,R).
$$
Since up to an $O(1)$ factor $\alpha_2(s)\,\alpha_{+}(s,|w|)$ is increasing in $s$,
the product of the first two terms is at most
$O(1)\,\alpha_2(|w|)$ and since $|w|\ge 2\,s\ge r$, the whole product is at most 
$$
O(1)\,\alpha_2(r)\,\alpha_{W_\theta}(r,R).
$$

\medskip\noindent
{Case 3(c)}: $|w| \ge R/2$; $s\le |w|/2$.
For $H$ to be visited, we need our 2 arms event holding 
within distance $s/2$ of $H$ and if 
$2s<d(p_0,w)$ it is also necessary that
a white crossing 
occurs between distances $2s$ and $d(p_0,w)\wedge |w|$ from
$w$.
(Note that the latter event takes place in the upper half plane.)
These are independent and since $|w|\ge R/2$ we get 
$$
\Pb{H\text{ visited }\md p_0}\le O(1)\,\alpha^+\bl(2s,d(p_0,w)\wedge (R/2)\br)\,\alpha_2(s/2)\,.
$$
As in Case 2, this is compatible with~\eref{e.revealment}.

\medskip
This covers all possible cases, and hence establishes that the
revealment is as claimed.

\medskip

We now proceed to discuss the algorithm and the revealment when $k>1$.
It turns out simplest in fact to
modify the event $A_{W_\theta}^k(r,R)$ as follows. Partition  the
outer boundary $\p_R D$ into $k$ arcs  of roughly equal diameter
$Y_1,Y_2,\ldots,Y_k$ (ordered counterclockwise) and let $\tilde{A}_{W_\theta}^k(r,R)$
be the event that for every odd [respectively, even] $i\in\{1,2,\dots,k\}$ there is a black
[respectively, white]
crossing in $D$ from $\p_r D$ to $Y_i$.
Thus, instead of looking at the set of times for which $A_{W_\theta}^k(r_0,R)$
occurs (where $r_0=10\,k$, say), we will look at the set of times at which $\tilde A_{W_\theta}^k(r_0,R)$
occurs.
Clearly, $\tilde{A}_{W_\theta}^k(r,R)\subset A_{W_\theta}^k(r,R)$,
and therefore this is justified.
We will also use the relation
\begin{equation}\label{e.AA}
\Pb{{A}_{W_\theta}^k(r,R)} \le C^\theta_k\, \Pb{ \tilde A_{W_\theta}^k(r,R)}\,,
\end{equation}
for some constant $C_k^\theta$, depending only on $k$ and $\theta$,
which holds by Remark~\ref{r.primus}.

If $Y\subset \p_R D$ is an arc, let $A_Y^1(r,R)$
[respectively, $A_Y^{-1}(r,R)$] be the
event that there is a white [repectively, black] crossing from $Y$ to $\p_r D$
in $D$.
Suppose that each for each $i=1,2,\dots,k$, we have
 a partition $Y_i=Y_{i,+}\cup Y_{i,-}$
of $Y_i$ into two arcs $Y_{i,+}$ and $Y_{i,-}$.
Since $A_{Y_i}^{\pm 1}(r,R)=A_{Y_{i,+}}^{\pm 1}(r,R)\cup A_{Y_{i,-}}^{\pm 1}(r,R)$,
we have
\begin{equation}\label{e.union}
\tilde{A}_{W_\theta}^k(r,R)
=\bigcup_{y\in \{-,+\}^k}
\bigcap_{i=1}^k A_{Y_{i,y_i}}^{{(-1)}^i}(r,R)\,.
\end{equation}

The algorithm starts out by picking points $x_i\in Y_i$,
randomly, uniformly and independently.
Then $Y_{i,+}$ and $Y_{i,-}$ are chosen as the two components
of $Y_i\setminus \{x_i\}$.
For each of the $2^k$ possible $y\in\{-,+\}^k$,
the algorithm then proceeds to determine if
the corresponding component
$$
A(y):=
\bigcap_{i=1}^k A_{Y_{i,y_i}}^{{(-1)}^i}(r,R)
$$
of~\eref{e.union} has occured.
For that purpose, interfaces are started at each of the points $x_i$,
and are followed until the event has been determined one way or the other.
(Of course, the interface will have either white on the left and
black on the right or vice versa, depending on the color of crossing
it is meant to detect and whether the corresponding arc
$Y_{i,\pm}$ is to the left or right of $x_i$.)
However, the order in which the interfaces are extended is somewhat important.
A simple rule that works is that among the hexagons that are necessary
to extend the $k$ interfaces one more step, the algorithm chooses the one that is
farthest away from $0$.
The event $A(y)$ is decided positively only if
all $k$ interfaces reach $\p_r D$.

The revealment of this algorithm is at most $k\,2^k$ times the maximum
probability that the interface started at $x_i$ visits a hexagon
$H$ before the determination of the corresponding
$A(y)$ is terminated. Here, the maximum is over all hexagons $H\subset D$
and all $i\in\{1,2,\dots,k\}$.
The corresponding bound is attained as in the case
$k=1$, but now $\alpha_{W_\theta}^1$ is replaced by
$\alpha_{W_\theta}^k$. 
Our rule of thumb for selecting which interface to
extend guarantees that we never examine a hexagon
$H$ unless $\tilde A_{W_\theta}^k\bl(\dist(0,H)+O(1),R\br)$ occurred.
As in the case $k=1$, when estimating the revealment it is important that
$\alpha_2(r,R)\le O(1)\,\alpha_{W_\theta}^k(r,R)$.
In the range $\theta> 4\,\pi\,k\,(k+1)/3$, which is the
relevant range for the lower bound in Theorem~\ref{th:wedges}, this follows
from~\eref{e.awt}.

The remainder of the proof goes through as before.
\QED

\proofof{lower bound in Theorem~\ref{th:cones}}
Here we simply say that the proof for the lower bounds in 
Theorem~\ref{th:wedges} can be carried out in a similar way. 
In fact, for $k\ge 2$, the proof
is simpler topologically than the $k=1$ case for the plane,
since we do not need to worry about interfaces making complete circuits
around the origin 
(if this ever happens, the event in question cannot occur and we stop the 
algorithm).
\QED

\section{Upper bounds for $k$-arm times} \label{sec:wedgesupper}

The following result,
which will be useful for the proofs
of the upper bounds in Theorems~\ref{th:wedges} and~\ref{th:cones},
is abstract: the graph structure does not play any role.
Let $A$ be an event involving independent Bernoulli $(1/2,1/2)$ random variables
$X_1,X_2,\ldots,X_m$. Recall that the influence of the index $i$ on $A$, denoted $I_i(A)$,
is the probability that $X_i$ is pivotal; namely,
that changing the value of $X_i$ changes whether $A$ 
occurs or not. 
The sum of the influences is denoted by $I(A)=\sum_i I_i(A)$.

\begin{theorem} \label{th.hdgeneral}
Let $\{A_n\}_{n\ge 1}$ be some sequence of events 
in $\{0,1\}^V$, each depending on only finitely many 
coordinates. 
Assume that $\lim_{n\to\infty} \Ps{A_n} = 0$.
Let $\omega_t$ be the Markov process on $\{0,1\}^V$ where
independently $0$'s go to $1$ at rate $1/2$,
$1$'s go to $0$ at rate $1/2$ and started according
to its stationary distribution $\pi_{\frac{1}{2}}$.
Let $T$ be the set of exceptional times $t$
at which $\omega_t\in\bigcap_{n\ge 1}A_n$.
If $\liminf_{n\to\infty} I(A_n) < \infty$, then 
$T=\emptyset$ a.s. Otherwise, the Hausdorff dimension of $T$ is a.s.\ at most
\begin{equation}\label{e.expression}
\liminf_{n\to\infty}
\left({{1-\frac{\log \Ps{A_n}}{\log I(A_n)}}}\right)^{-1}\,.
\end{equation}
\end{theorem}

\proof
Let $T_n:=\{t\in [0,1]:\omega_t\in A_n\}$, let $\partial T_n$ be the boundary
points of $T_n$ in $(0,1)$ and set $N_n:=|\partial T_n|$. We claim that
\begin{equation}\label{e.Nn=In} 
\Eb{N_n}={I(A_n)}/{2}.
\end{equation}
To see this, write
$N_n=\sum_v N^v_n$ where $N^v_n$ counts the number of elements in
$\partial T_n$ at which time the vertex $v$ flipped. We now need to
show that, for each vertex $v$, $\Eb{N^v_n}$ is $I_v(A_n)/2$.
Given a time interval $[t,t+dt]$, the probability that 
there is a time point in the interval which contributes to $N^v_n$
is equal to $I_v(A_n)\,dt/2 + o(dt)$ and the probability of $k\ge 2$ such time
points is clearly $O(dt^k)$. From this,~(\ref{e.Nn=In}) easily follows.

For any $\eps > 0$, let $T^\eps_n$ be the $\eps$-neighborhood of
$T_n$ intersected with $[0,1]$.
Since $T^\eps_n\subseteq T_n\cup \bigcup_{x\in \p T_n}[x-\eps,x+\eps]$,
\begin{equation}\label{e.star} 
\mu(T^\eps_n)\le \mu(T_n)+ 2\,N_n\,\eps\,,
\end{equation}
where $\mu$ denotes Lebesgue measure.
For any set $U$ and $\epsilon>0$, let $\calN(U,\epsilon)$ denote the number of 
$\epsilon$ intervals needed to cover $U$.
From the above, using the fact that the intervals comprising
$T^\eps_n$ all have length at least $\eps$, it follows that
$\calN(T^\eps_n,\eps)\le2\,\mu(T^\eps_n)\,\eps^{-1}$, and so,
using~\ref{e.star},
$$
\calN(T_n,\eps)\le \calN(T^\eps_n,\eps)\le
{2\,\mu(T_n)\,}{\eps^{-1}}+ 4\,N_n\,.
$$
Therefore, by Fubini's theorem and~\ref{e.Nn=In}, 
\begin{equation}\label{e.covernumber}
\Eb{\calN(T_n,\eps)}\le {2\,\Ps{A_n}}\,{\eps^{-1}}+ 2\,I(A_n)\,.
\end{equation}

We now temporarily assume that $\liminf_{n\to\infty} I(A_n) =\infty$. 
Let $a_n=\Pb{A_n}/I(A_n)$, which goes to $0$ as $n\to\infty$.
By (\ref{e.covernumber}), we have 
\begin{equation}\label{e.neededbound}
  \Eb{\calN(T,a_n)}\le 
  \Eb{\calN(T_n,a_n)}\le 4\,I(A_n)\,.
\end{equation}
By passing to a subsequence if necessary, we assume with no loss of generality
that the $\liminf$ in~\eref{e.expression} is a limit. Let $L$ denote the value
of that limit.
It is elementary to check that for every $\eps>0$, for all sufficiently large $n$,
$$
I(A_{n})\le \left(\frac{I(A_{n})}{\Ps{A_{n}}}\right)^{L+\epsilon}.
$$
This together with (\ref{e.neededbound})
implies that  %
the Hausdorff dimension of $T$ is at most $L+\epsilon$ a.s. As 
$\epsilon$ is arbitrary, this completes the proof in the case
$I(A_n)\to\infty$.

Since $T_n\ne \emptyset$ implies that $N_n\ge 1$ or $T_n\supseteq (0,1)$,
it follows by (\ref{e.Nn=In}) and Markov's inequality that 
$$
\Pb{T_{n} \neq \emptyset}\le \Pb{A_{n}} + I(A_{n})\,.
$$
Thus, $T=\emptyset$ a.s.\ when $\liminf_n I(A_n)= 0$.

The case $\liminf_n I(A_n)\in (0,\infty)$ requires a different argument.
Let $\eps_n=\sqrt{\Ps{A_n}}$. By (\ref{e.covernumber}), we have 
$
\liminf_{n\to \infty}\Eb{\calN(T_n,\eps_n)} < \infty
$.
Since $\eps_n\to 0$, the cardinality $|T|$ of $T$ is bounded from above by
$ \liminf_{n\to\infty}\calN(T_n,\eps_n) $.
Fatou's lemma yields $\Eb{|T|}< \infty$ and
hence $\Pb{|T|<\infty}=1$. We finally conclude that
$\Ps{T\neq\emptyset}=0$ by combining
\cite[Theorem 6.7]{GS} and \cite[(2.9)]{FG}.
\QED

\proofof{Theorem~\ref{th:cones}}
Since the lower bounds have been established in Section~\ref{sec:wedgeslower},
it remains to prove the upper bounds.
Fix $k=1$ or $k>1$ even.
Let $A_R$ be the event that there are $k$ disjoint
crossings of the annulus $D_R:=\{z\in C_\theta: 10k\le |z|\le R\}$, where we require
that the colors be alternating if $k\ne 1$.
Here, $|z|$ denotes the distance to $0$, which is the apex of the cone $C_\theta$.
One can prove that 
\begin{equation}\label{e.coneexpo}
\Pb{A_R}
=\begin{cases}
R^{-5\pi/(24\theta)+o(1)}&k=1,\\
R^{{(1-k^2)\pi}/({6\theta}) + o(1)}& k>1,
\end{cases}
\end{equation}
in the very same way that we have justified~\eref{e.awt}.
By Theorem~\ref{th.hdgeneral} (and easy algebraic manipulation), it therefore suffices to prove
that
\begin{equation}\label{e.pivotal}
I(A_R)/\Ps{A_R}\le R^{3/4+o(1)}.
\end{equation}
Let $H$ be a hexagon in $C_\theta$, and let $s=s(H)$ be the distance from
$H$ to $0$.
For $H$ to be pivotal 
it is necessary that there would be
$k$ disjoint (alternating, if $k>1$)
crossings from distance $10\,k$ to
$s/2$ from the origin (unless $s/2\le 10\,k$)
and between distances $2\,s$ and
$R$ (unless $2\,s\ge R$).
Likewise, there should be $4$ alternating crossings
between $H$ and distance
$(s/2)\wedge\dist(H,\p D_R)$ from $H$.
These events are independent.
Using the quasi-multiplicative property
of the $k$-arm crossing events (Remark~\ref{r.qmc})
and~\eref{e.fullpk} with $k=4$,
this gives when $s<8\,R/9$
\begin{equation}\label{e.sumI}
I_H(A_R)\le O(1)\,\Ps{A_R}\,s^{-5/4+o(1)}.
\end{equation}
where this $O(1)$ factor (as well as those appearing
below) may depend on $k$ and $\theta$.
Since the number of hexagons in $C_\theta$
satisfying $s=s(H)<\rho$ is $O(\rho^2)$,
an easy calculation yields
$$
\sum_{H:s(H)< 8R/9} I_H(A_R)\le O(1)\,\Ps{A_R}\,R^{3/4+o(1)}\,.
$$

Now suppose that $H$ is a hexagon satisfying
$s(H)\ge 8\,R/9$.
For $H$ to be pivotal, it is necessary that there would
be $k$ (alternating, if $k>1$) crossings in $C_\theta$
between $\{|z|=10\,k\}$ and $\{|z|=R/2\}$,
there should be $4$ alternating crossings
between $H$ and distance $\dist(H,\p D_R)/2$ from $H$,
and there should be $3$ alternating
crossings
between distance $2\,\dist(H,\p D_R)$ and
distance $R/2$ from a point on $\p D_R$ closest to $H$.
The latter event is governed by the $3$-arm half
plane exponent, whose asymptotic behaviour is described by~\eref{e.halfpk}.
Since $s+\dist(H,\p D_R)=R+O(1)$, we get
$$
I_H(A_R)\le O(1)\,\Ps{A_R}\,(R-s)^{-5/4+o(1)}\,((R-s)/R)^{2+o(1)}\,.
$$
Since for $b\ge 1$ there are $O(b\,R)$ hexagons at distance $\le b$ from
$\{|z|=R\}$, another easy calculation gives
$$
\sum_{H:s(H)\ge  8R/9} I_H(A_R)\le O(1)\,\Ps{A_R}\,R^{3/4+o(1)}\,.
$$
Together, this yields~\eref{e.pivotal} and the proof is complete.
\QED

\proofof{Theorem~\ref{th:wedges}}
The lower bound was proved in Section~\ref{sec:wedgeslower},
and so only the upper bound needs to be justified.
The proof proceeds like the proof of the upper bound in Theorem~\ref{th:cones},
except that the influence estimates are slightly different.

Let $D_R=\{z\in W_\theta:10\,k\le |z|\le R\}$,
$A_R$ be the $k$-arm event in $W_\theta$
between $\{z:|z|=10\,k\}$ and
$\{z:|z|=R\}$, and $H\subset D_R$ be a hexagon.
Let $s=s(H)=\dist (0,H)$, and let $b=b(H)=\dist (H,\p D_R)$,
where we write $\p D_R$ for the boundary of $D_R$ in $C_\infty$,
i.e., the points on $\p W_\theta$ are included.
For $H$ to be pivotal for $A_R$, it is necessary
that the $k$ arm event holds
between distance $10\,k$ and $s/2$ from $0$ (unless
$s/2\le 10\,k$) as well as between
distances $2\,s$ and $R$ (unless $2\,s\ge R$),
that the alternating $4$-arm
event hold between $H$ and distance $b/2$ away from $H$,
and that the alternating $3$-arm event must hold
between distances $2\,b$ and $s/4$ away from a point
in $\p D_R$ closest to $H$ (unless $2\,b\ge s/4$).
There are $O(b's')$ hexagons $H$ satisfying $b(H)\le b'$
and $s(H)\le s'$.
The rest of the proof proceeds like that of Theorem~\ref{th:cones},
and is left to the reader.
\QED

\proofof{Theorems~\ref{th.no2clusters},
\ref{th.differentclusters},
\ref{th.halfplane} and
\ref{th.halfplane2}}
At any time at which there are 2 infinite white clusters in the plane,
 we must also have
the 4-arm event occuring (with alternating colors) but by
Theorem~\ref{th:cones} (with $k=4$ and $\theta=2\pi$), there are no such times.
This proves Theorem~\ref{th.no2clusters}.

At any time at which there are 2 infinite different colored clusters, we must also 
have the 2-arm event occuring (with different colors) but by
Theorem \ref{th:cones} (with $k=2$ and $\theta=2\pi$), the set of
such times has Hausdorff dimension at most $2/3$.
This proves Theorem~\ref{th.differentclusters}.
The other two theorems are similarly proved.
\QED

\section{The square lattice} \label{sec:square}

We start this section by proving 
Theorem~\ref{th.no3clustersZ2}.
Afterwards, possible ways in which our arguments for
Theorem~\ref{th:except} may be improved to apply to $\Z^2$ 
as well, will be discussed.

In the proof of Theorem~\ref{th.no3clustersZ2} we will
use the fact that the $6$ alternating arms exponent is
larger than $2$, or, more precisely, 
that the probability for $6$ alternating arms
between radii $r$ and $R$ is
bounded above by $O(1)\,(r/R)^{2+\eps}$ for some $\eps >0$.
This is essentially due to~\cite[Lemma 5]{KSZ},
but a proof is also given in the appendix (Corollary~\ref{c.56}).

\proofof{Theorem~\ref{th.no3clustersZ2}}
For $0<r<R$, let $S(r,R)$ be the event
that there are $3$ different clusters that connect the circles
of radii $r$ and $R$ about $0$.
By the above mentioned bound on the alternating $6$-arm probabilities,
We may choose some fixed $\eps>0$ and
some function $\rho=\rho(r)>r$ such that for static critical bond percolation on $\Z^2$,
for all $r$,
\begin{equation}\label{e.rhochoice}
\Pb{ S(r,\rho)}\le\rho^{-2-\eps}\,.
\end{equation}
Consider some bond $e$, and let $F(e)$ be the event that
$e$ is pivotal for $S(r,\rho)$.
Then $\Pb{F(e)}$ is just the influence $I_e(S(r,\rho))$.
Assume that $\Pb{F(e)}\ne 0$.
Note that the events $F(e)$ and $\{e\text{ is open}\}$ are independent
events. This implies that $\Pb{S(r,\rho)\md F(e)}=1/2$,
which one may write
$ \Pb{F(e)\cap S(r,\rho)}=\Pb{F(e)\cap\neg S(r,\rho)}$.
Since this applies to every bond $e$, we conclude that
the expected number of pivotals on the event $S(r,\rho)$
is half of the total influence $I(S(r,\rho))$.
However the number of pivotals for $S(r,\rho)$ is bounded
by the total number of edges intersecting the
disk of radius $\rho$ about the origin, which is
certainly $O(\rho^2)$. Thus, 
$$
I(S(r,\rho))\le 2\,\Pb{S(r,\rho)}\,O(\rho^2) =O(\rho^{-\eps}).
$$
Consequently, by Theorem~\ref{th.hdgeneral}, for every $r_0>0$ a.s.\ there
are no exceptional times in which
$\bigcap_{r>r_0} S(r,\rho(r))$ holds.
This proves our theorem.
\QED

\begin{remark}
An alternative way to prove the above result is based on using the fact that the 6-arm 
exponent is strictly larger than 2 together with the fact that the number of different
configurations (counting repetitions) that appear in a ball of radius 
$n$ during the time interval $[0,1]$ has a Poisson distribution
with a parameter which is at most $O(1) n^2$.
\end{remark}

\bigskip

As we will briefly explain below, the proof of
Theorem~\ref{th:except} almost works for bond percolation on the square
grid. In fact, there are several alternative routes by which
the result might perhaps be extended to $\Z^2$:
\begin{enumerate}
\item establishing 
\begin{equation}\label{e.exponeed}
\alpha_2(r)\le O(r^{-\eps})\,\alpha(r)^{2}
\end{equation}
 for $\Z^2$ for some fixed $\eps >0$,
\item improving the estimate~\eref{e.noise},
\item proving the existence of an algorithm
(or a witness which would still permit the use of
Theorem~\ref{t.noise}) with smaller revealment,
\item extending Smirnov's theorem to $\Z^2$.
\end{enumerate}

Note that the weaker version of~\eref{e.exponeed}
$\alpha_2(r)\le \alpha(r)^2$ follows from either the Harris-FKG 
inequality or Reimer's inequality~\cite{Reimer}.
Kesten and Zhang have proved some related strict inequalities
between exponents~\cite{KZ}, but it seems that their
methods are not sufficient to prove~\eref{e.exponeed}.

We now explain why~\eref{e.exponeed} in the $\Z^2$ setting implies
exceptional times for $\Z^2$.
First we want to have the revealment
for the algorithm determining $\croRr$ bounded by
$O(1)\,\alpha_2(r)\,\alpha(r,R)$.
One problem seems to be that
the bound on the revealment
for the triangular grid involves the
summand featuring $\alpha^+$, which is relatively negligible, while
on $\Z^2$, we do not know how to prove that the other summand
dominates. The fix is to replace the deterministic $R$
by a random $R_t'\in[R,2R]$.
The random variable $R_t'$ will depend on some extra random bits,
that we add, and these random bits also evolve in time.
We construct the dependence of $R'$ on these bits
so that $R'$ can be calculated by an algorithm with very
small revealment.
This is rather easy to arrange, because we are not limited 
in the number of bits that we may take. 
If we consider an edge whose distance from the origin $a$ is
in the range $[R/2,2R]$, then the probability
that the edge is examined given $R'$ is at most
$O(1)\,\alpha_2(R'-a)\,1_{R'\ge a-1}$.
By~\eref{e.qm}, this is at most
$O(1)\,\alpha_2(R)\,\alpha_2(R'-a,R)^{-1}\,1_{R'\ge a-1}$.
The probability that $|R'-a|\le 2^j$ is at most $O(1)\,2^{j}/R$.
We also know that $\alpha_2(r_1,r_2)^{-1}\le O(1)\,(r_2/r_1)^{1-\eps'}$
for some $\eps'>0$, by Reimer's inequality~\cite{Reimer} and~\eref{e.fivearm}.
It follows that the probability that such an edge is examined is 
$O(1)\,\alpha_2(R)$.
The $r^{o(1)}$ factor in Theorem~\ref{t.annulus}
is easily replaced by an $O(1)$ factor, if we use Proposition~\ref{p.qm}
in the course of the proof.
Then we get~\eref{e.zz} for the square grid, but without the $r^{o(1)}$ factor.
We may then choose the dependence
between $r$ and $t$ such that
$\alpha(r)\approx t\,r^{\eps/2}$, where $\eps$ is the constant
in~\eref{e.exponeed}. The rest is immediate from~\eref{e.zz}, since
clearly $r^{-\eps/2}\le O(1)\,t^{\eps'}$ for some $\eps'>0$.

A consequence of this argument, which applies without assuming~\eref{e.exponeed},
is that for bond percolation on $\Z^2$ we have
$$
\Pb{V_{t,R}\cap V_{0,R}} \le O(t^{-1})\,\Pb{V_{t,R}}^2.
$$
This gives yet another illustration as to how close the result for
$\Z^2$ seems to be --- if the $t^{-1}$ term was improved to
$t^{-1+\eps}$, that would have been enough.
Consequently, significant improvements in the algorithm or in~\eref{e.noise}
would also be sufficient.

\section{Some open questions} \label{sec:questions}

Following are a few questions and open problems suggested by the
present paper.

\begin{enumerate}
\item For the results in Theorems \ref{th:wedges} and \ref{th:cones},
what is the Hausdorff dimension of the set of exceptional times 
in question? We tend to believe that the answer is the upper bound.
In particular, is the upper bound of $31/36$ in
Theorem \ref{th.hd} the correct answer?
\item Prove that there exist exceptional times for percolation
on the square lattice. (See Section~\ref{sec:square} for a discussion.)
\item What is the best $\gamma$ for which Theorems~\ref{th:crossingquantsquare}
 and~\ref{th:crossingquant} hold?
\item
What is the best revealment of an algorithm determining
the event $\cro_R$ in Theorem~\ref{t.square}?
\item
What is the sharp form of Theorem~\ref{t.noise}?
\item
What are the properties of the infinite cluster at an exceptional
time at which it exists? 
For example, what is the growth rate of the number of vertices
in the Euclidean disk of radius $r$ around the origin which
belong to the cluster of the origin at the first time
$t\ge 0$ in which the cluster is infinite?
Is the growth rate the same at all exceptional times?
\item
What is the relationship between the exceptional infinite
cluster and the incipient infinite cluster?
\end{enumerate}

\appendix
\section{Appendix: Quasi-multiplicativity}\label{s.appendix}

In this appendix, we discuss the $k$-arm probabilities
and prove that they satisfy the corresponding analogue of 
the relation~\eref{e.multip}.
For $R>0$, let $H_R$ be the union  of the hexagons intersecting $B(0,R)$.
For $R>r>0$ let $A_j(r,R)$ denote the event that
there are at least $j$ crossings from $\p H_r$ to $\p H_R$, of alternating colors.
The following result refers to critical site percolation on the triangular
grid and critical bond percolation on the square grid.

\begin{proposition} \label{p.qm}
Let $j>0$ be even. There is a constant $C$, depending only on $j$,
such that for all $r<r'<r''$
\begin{equation}\label{e.qm}
C^{-1}\, \Pb{A_j(r,r'')}\le
\Pb{A_j(r,r')}\,\Pb{A_j(r',r'')}\le C\, \Pb{A_j(r,r'')},
\end{equation}
and $\Pb{A_j(r,2\,r)}>1/C$ if $\Pb{A_j(r,2\,r)}>0$
(i.e., if $r$ is large enough to allow $j$ crossings).
Moreover, a corresponding statement holds for critical bond percolation on the
square grid which alternate between primal and dual crossings.
\end{proposition}

This theorem would have been a useful tool in~\cite{SW},
had it been available. Instead, the authors of that paper 
proved a weaker form of this which was good enough for their
purposes. Our proof below uses techniques from~\cite{K2},\cite{LSWup2}
and \cite{SW}.
Indeed, the entire results of the appendix do follow from the ideas
of~\cite{K2}. We include them here for completeness, and for ease
of reference. Additionally, though the basic ideas are the same,
in several respects our treatment is a bit different from~\cite{K2}.

Below, we will work in the setting of the triangular grid.
The proof for the square grid is essentially the same.
In the setting of the triangular grid, there is the color exchange trick~\cite{KSZ,ADA},
which shows that the probability for having alternating crossings is comparable to the
probability of any color sequence as long as both colors are present.
In the setting of the square grid, as far as we know, such a trick does not exist.
At the end of the appendix we will explain how the proof of Proposition~\ref{p.qm}
can be generalized to any color sequence.

In the following, an interface
from $\p H_r$ to $\p H_R$ is an oriented simple path in the
hexagonal grid that has one color of hexagons adjacent to it on the
right, and the opposite color adjacent to it on the left.
Thus, it is the common boundary of a black crossing and a white crossing.

For $R>r>1$, consider all the interfaces crossing from $\p H_{r}$ to
$\p H_{R}$, and define $s(r,R)$ to be the least distance between any
pair of endpoints of two interfaces on $\p H_R$.
If there are no interfaces, we take $s(r,R)=\infty$.
Note that $s(r,R)$ is monotone non-increasing in $r$.
This quantity will roughly measure the \lq\lq quality\rq\rq\
of the interfaces; when $s(r,R)$ is comparable to $R$, the interfaces
are well separated, and, as we will see, easier to extend.

\begin{lemma}\label{l.apart}
For all $a\in(0,1)$, $R>0$, $\delta>0$
$$
\Pb{s(a\,R,R)<\delta\,R}\le C\,\delta^\eps,
$$
where $C=C(a)$ is a constant depending only on
$a$ and $\eps>0$ is an absolute constant.
\end{lemma}

The lemma probably follows from \cite[Lemma 2]{K2},
but since the proof is rather short, we include a proof
for completeness.

\proof
We prove this in the case $a=1/2$.
The general case is essentially the same.
Let $\alpha\subset\p H_R$ be an arc of diameter 
$R/3$, and let $Y$ be the set of points
in $H_R$ at distance at most $R/3$ from
$\alpha$.
Let $\alpha_1$ be one of the two
arcs in $\p Y\cap\p H_R\setminus\alpha$.
Let $\beta_1,\beta_2,\dots,\beta_k$
be the interfaces crossing
from $\p Y\setminus\p H_R$ to $\alpha$,
ordered so that for $i_1<i_2\le k$,
the interface $\beta_{i_1}$ separates
$\alpha_1$ from $\beta_{i_2}$ in $\closure Y$.
Fix a positive integer $i$
and condition on $i\le k$ and on $\beta_i$.
Let $\hat\beta_i$ denote the union of the hexagons adjacent
with $\beta_i$.
Then the
percolation in the connected component $Y_i$ of
$Y\setminus\hat \beta_i$ separated from $\alpha_1$
by $\beta_i$ remains unbiased.
Suppose that the hexagons in $\hat\beta_i$
adjacent to $Y_i$ are white, say.
Let $z_i$ denote the endpoint of $\beta_i$ on
$\alpha$.
The RSW theorem implies that
there is some constant $\eps>0$ such that with conditioned
probability $1-O(1)\,\delta^\eps$
there is a white crossing 
in $Y_i\setminus B(z_i,\delta\,R)$ from
$\p\hat\beta_i$ to $\p H_R$.
On that event, it is clear that
if $k\ge i+1$, then $|z_i-z_{i+1}|\ge \delta\,R$.
We conclude that
$$
\PB{k\ge i+1,\,|z_i-z_{i+1}|\le\delta\,R\md
k\ge i,\,\beta_i}\le O(1)\,\delta^\eps.
$$
The RSW theorem also implies that there is conditioned probability
bounded away from zero that $k=i$
given $k\ge i$ and $\beta_i$ (this would be guaranteed by
an appropriate crossing in $Y_i$ from $\p\hat\beta_i$
to $\p H_R\setminus\alpha$). Therefore,
$\Pb{k\ge i}\le c^i$ for some constant $c<1$.
Thus, 
\begin{multline*}
\PB{k\ge i+1,\,|z_i-z_{i+1}|\le\delta\,R}
\\
=
\PB{k\ge i+1,\,|z_i-z_{i+1}|\le\delta\,R\md k\ge i}\,
\Pb{k\ge i}= O(1)\,c^i\,\delta^\eps\,.
\end{multline*}
We sum this over all $i=1,2,\dots$, and over an appropriate
covering of $\p H_R$ by $O(1)$ arcs $\alpha$ of diameter $R/3$.
The lemma follows.
\QED

Next, we prove a statement
saying, roughly, that if the crossings are \lq\lq reasonably good\rq\rq,
then there is a conditioned probability bounded away
from zero that they extend and
the extensions are \lq\lq very good\rq\rq.

\begin{lemma}\label{l.improve}
For every $j>0$ even, there is a constant $\bar\delta=\bar\delta(j)>0$
such that 
for every $\delta>0$ there is some constant $c(\delta)>0$, depending only on $\delta$,
such that when $R>r$,
$$
\PB{A_j(r,4\,R)\cap\{s(r,4R)>4\,\bar\delta\,R\}\md  A_j(r,R)\cap\{ s(r,R)>\delta\,R\} } > c(\delta)\,.
$$
\end{lemma}

\proof
Set $S=\closure{H_R\setminus H_r}$.
We assume that $A_j(r,R)$ holds and
that $s(r,R)>\delta\,R$.
Let $\gamma_0,\dots,\gamma_{k-1}$  ($k\ge j$) be the collection of all interfaces crossing
from $\p H_r$ to $\p H_R$, in counterclockwise order,
where we choose the indexes so that $\gamma_0$ has white hexagons
on the right hand side. 
(The interfaces are oriented from $\p H_r$ to $\p H_R$.)
In the following, we set $\gamma_i:=\gamma_{i'}$
when $i\notin \{0,1,\dots,k-1\}$ and $i'=i\mod k$.
Set $\Gamma=\bigcup_{i\in \N}\gamma_i$.
How does conditioning on the interfaces $\gamma_0,\dots,\gamma_{k-1}$
affect the percolation process? 
Note that the fact that there are no more than $k$
interfaces means that for each $i\in\N$ there is a white crossing
in $S\setminus\Gamma$
from the right hand side of $\gamma_{2i}$ to
the left side of $\gamma_{2i-1}$ and 
a black crossing in $S\setminus\Gamma$ from the left side of
$\gamma_{2i}$ to the right side of $\gamma_{2i+1}$.
Otherwise, the configuration is unbiased on hexagons
that are not adjacent to these interfaces.

\begin{figure}
\SetLabels
\E(.5*.5)$H_R$\\
\T(.5*.98)$\p H_{4R}$\\
\endSetLabels
\centerline{\epsfysize=2.5in%
\AffixLabels{%
\epsfbox{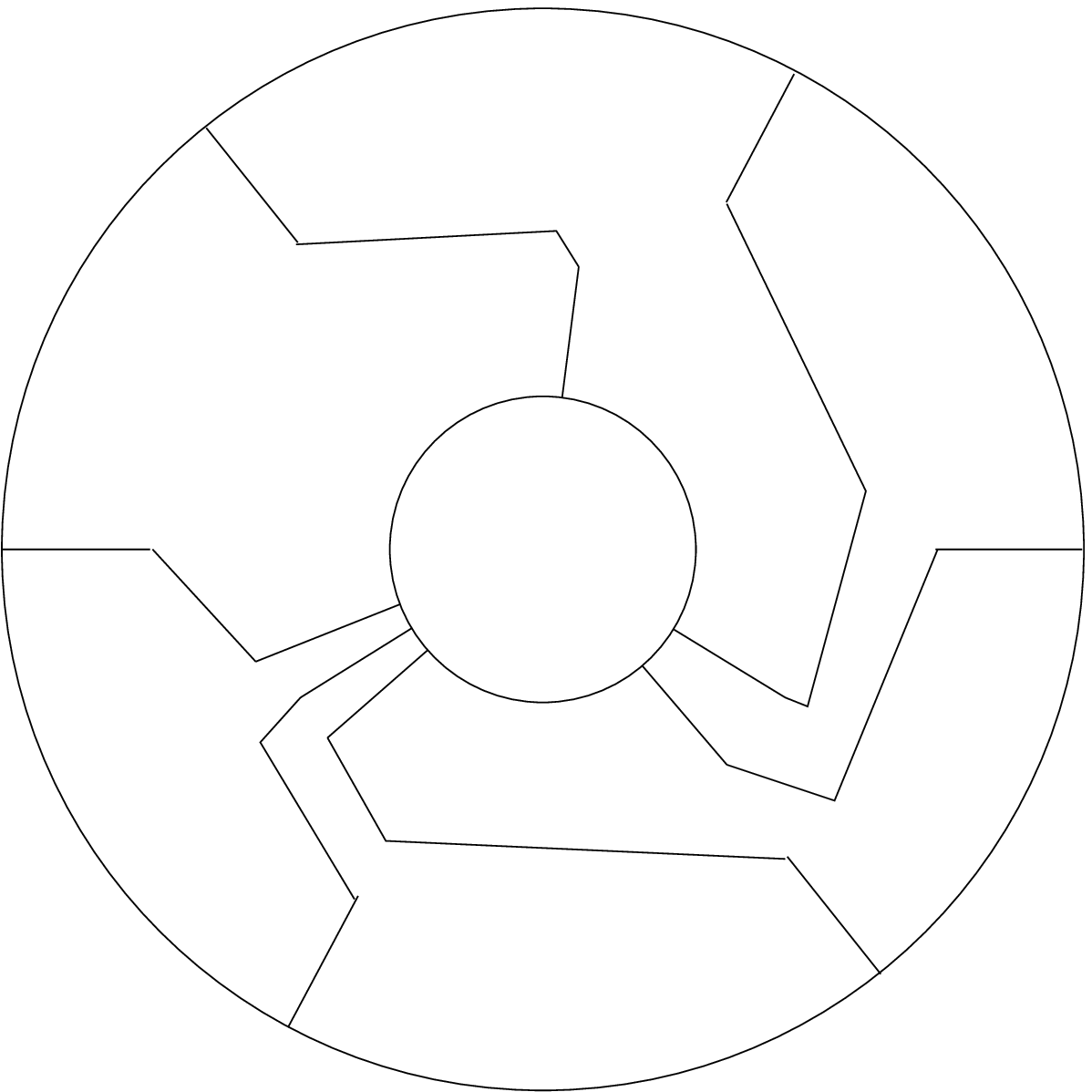}%
}%
}
\begin{caption} {\label{f.betas}The paths $\beta_i$.}
\end{caption}
\end{figure}

Let $z_i$ be the endpoint of $\gamma_i$ on $\p H_R$, $i\in\N$.
For $i=0,1,\dots,{j-1}$, let $w_i$ be a point $\p H_R$ that
is roughly in the center of the counterclockwise arc from
$z_i$ to $z_{i+1}$. Then $|w_i-z_i|\ge\delta\,R/5$,
 and the same is true for $|w_i-z_{i+1}|$.
It is easy to see that there exist disjoint simple paths
$\beta_0,\dots,\beta_{j-1}$ satisfying the following.
(See Figure~\ref{f.betas}.)
(1) Each $\beta_i$ is a path in $\closure {H_{4R}\setminus H_R}$ 
from $w_i$ to a point $w_i'\in\p H_{4R}$.
(2) The points $w_i'$ are roughly equally spaced on
$\p H_{4R}$.
(3) Each $\beta_i\cap H_{2R}$ is contained in the line
through the origin containing $w_i$,
and each $\beta_i\setminus H_{3R}$ is contained
in the line through the origin containing $w_i'$.
(4) The distance from each of these paths to any other
path is at least $c_1\,\delta\,R$,
where $c_1\in(0,1/5)$ is some universal constant.
(5) Each $\beta_i$ has length at most constant
 times $R$, where the constant may depend on $j$.
For each $i\in\{0,1,\dots,j-1\}$ let
$\alpha_i$ be the arc of a circle whose center
is $z_i$, that has $w_i$ as endpoint,
that has the other endpoint on $\gamma_i\cup\gamma_{i+1}$,
that is otherwise disjoint from $\gamma_i\cup\gamma_{i+1}$
and is contained in $S$.

Let $\hat\beta_i$ be the connected component containing $w_i'$ of
the complement of $\Gamma$ in the $c_1\,\delta\,R/20$ neighborhood
of $\beta_i\cup\alpha_i$.
If $i\in\{0,1,\dots,j-1\}$ is odd, let $F_i$ denote the event that
there is a white crossing
in $\hat\beta_i$ from $\Gamma$ to $\p H_{4R}$.
Similarly, if $i$ is even, let $F_i$ denote the event that
there is a black crossing in $\hat\beta_i$ from
$\Gamma$ to $\p H_{4R}$.
It is easy to see that if $\bigcap_{i=0}^{j-1}F_i$ holds,
then $A_j(r,4R)$ holds as well.
The RSW theorem implies that $\Pb{F_i}$ is bounded from
below (depending on $\delta$) for a percolation process that is unbiased.
But, as we have seen,
the percolation on $S$ between $\gamma_i$ and $\gamma_{i+1}$
is only conditioned on having a crossing of the appropriate
color. By the Harris-FKG inequality, this is positively
correlated with $F_i$. By independence on
disjoint sets, the different $F_i$ are independent
given $\Gamma$
(assuming, as we may, that the distance between the different
sets $\hat\beta_i$ is significantly larger than the scale of the lattice). 
We conclude that for some $c(\delta)>0$,
$$
\PB{A_j(r,4\,R)\md  A_j(r,R)\cap\{ s(r,R)>\delta\,R\} } > c(\delta)\,.
$$

\begin{figure}
\SetLabels
\L(.83*.85)$\p H_{4R}$\\
\E(.01*.01)$\hat\beta_i$\\
\L\B(.39*.7)$\p H_{3.5\,R}$\\
\endSetLabels
\centerline{\epsfysize=2.5in%
\AffixLabels{%
\epsfbox{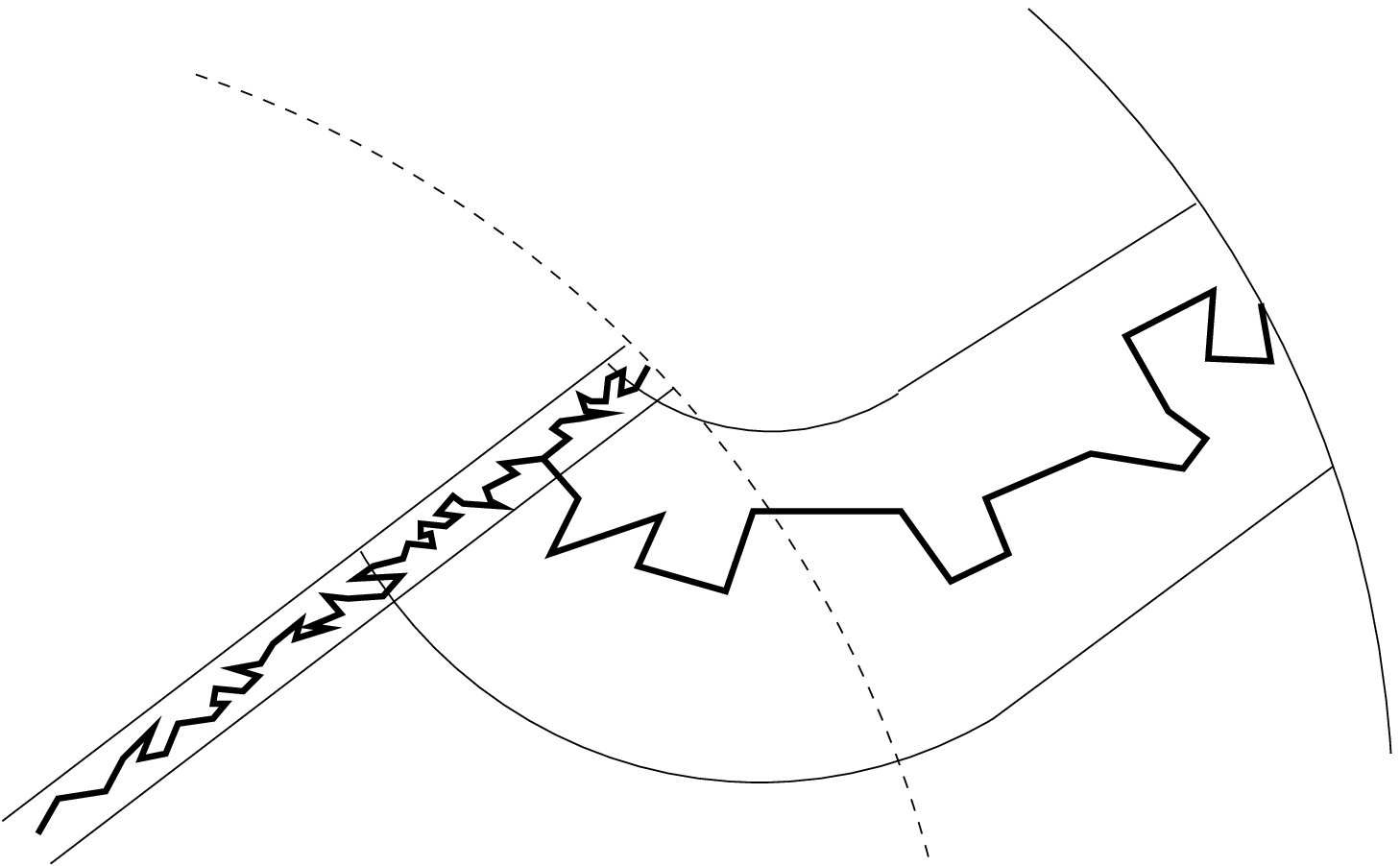}%
}%
}
\begin{caption} {\label{f.bent}A bent strip connecting with a channel.
Indicated are the leftmost crossing of
the intersection of the bent strip and the channel that connects
to $H_r$ and a reasonably likely crossing from it to $\p H_{4R}$ in the bent strip.}
\end{caption}
\end{figure}

Taking care of the condition $s(r,4R)\ge 4\,\bar\delta\,R$ is not too hard.
Suppose that in the above we truncate the paths
$\beta_i$ and the neighborhoods $\hat\beta_i$ by intersecting
them with $H_{3.5\,R}$.
We then condition on the \lq\lq leftmost\rq\rq\
crossing in $\hat\beta_i$.
See Figure~\ref{f.bent}. 
All this takes place within $H_{3.5\,R}$.
The conditional probability that these crossings in
the $\hat\beta_i$'s connect to
$\p H_{4\,R}$ is bounded away from
zero by a function of $j$ only (specifically, not $\delta$).
Thus, Lemma~\ref{l.apart} and the monotonicity
of $s(r,R)$ in $r$ shows that if
$\bar\delta=\bar\delta(j)>0$ is chosen small, with conditional
probability at least $1/2$ we are also likely to have $s(r,4R)\ge4\,\bar\delta\,R$,
as required.
\QED

Set 
$$
f(r,R):=\Pb{A_j(r,R)}, 
\qquad
g_\delta(r,R):=\Pb{A_j(r,R)\cap \{s(r,R)>\delta\,R\}}.
$$

\begin{lemma}\label{l.plenty}
There is a constant $\CCb(j)>0$, depending only on $j$,
such that for $R\ge 4\,r$
$$
\CCb (j)\, g_{\bar \delta}(r,R)\ge f(r,R)\,,
$$
where $\bar\delta=\bar\delta(j)$ is the constant introduced in
Lemma~\ref{l.improve}.
\end{lemma}

\proof
We assume that $f(r,R)>0$.
Let $\delta>0$ be small.
Set $N=\log_4(R/r)$ and let $m=m_\delta$ be the largest
integer in the range $0\le m\le N-1$
such that
$ g_\delta(r,4^{-i}R)\le f(r,4^{-i}R)/2 $
holds for every integer $i$ in the range
$0\le i< m$. 
Lemma~\ref{l.apart} and independence
on disjoint sets imply
$$f(r,R)-g_\delta(r,R)\le C\,\delta^\eps\, f(r,R/4),$$
and repeated applications of this inequality give
\begin{equation}\label{e.fR}
f(r,R)\le (2C)^m\,\delta^{\eps\, m} f\bl(r,4^{-m} R\br)\,.
\end{equation}
We claim that if $\delta$ is a sufficiently small positive constant, then
\begin{equation}\label{e.gd}
f(r,4^{-m}R) \le \CCc(j)\,g_\delta(r,4^{-m}R)
\end{equation}
for some constant $\CCc(j)$ depending only on $j$.
If $m\le N-2$, this follows with $\CCc(j)=2$ from the definition of $m$.
If $N-2<m\le N-1$, then RSW easily implies
$f(r,4^{-m}R)\ge \CCa (j)$ for some constant $\CCa(j)>0$,
and $f(r,4^{-m}R)-g_\delta(r,4^{-m}R)\le C\,\delta^\eps$ by Lemma~\ref{l.apart},
which gives~\eref{e.gd}.
On the other hand, repeated application of Lemma~\ref{l.improve}
gives
\begin{equation}\label{e.fullcircle}
c(\delta)\,c(\bar\delta)^{m-1}\,g_\delta(r,4^{-m}R)
\le g_{\bar\delta}(r,R) \le f(r,R)
\,.
\end{equation}
On combining this with~\eref{e.fR} and~\eref{e.gd}, one obtains
$$
c(\bar\delta)\,\CCc(j)\,\bl(2\,C\,\delta^\eps/ c(\bar\delta)\br)^m \ge c(\delta)\,.
$$
We choose $\delta$ sufficiently small so that $\delta^\eps< c(\bar\delta)/(4\,C)$.
Then the above shows that $m$ is bounded by a function of $\delta$ and $j$.
The proof is now completed by combining~\eref{e.fR},~\eref{e.gd}
and~\eref{e.fullcircle}.
\QED

\proofof{Proposition \ref{p.qm}}
The proof will be given only for the triangular grid,
since the proof in the setting of bond percolation on the
square grid is essentially the same.
As we remarked above, when $\Pb{A_j(r,2\,r)}>0$,
the RSW theorem easily gives $\Pb{A_j(r,2\,r)}>1/C(j)$,
for some $C(j)>0$.

We now assume that $r'>8\,r$ and $r''>8\,r'$.
Suppose that $A_j(r,r'/2)\cap\{s(r,r'/2)>\bar\delta\,r'/2\}$
holds. We also assume that the corresponding event
occurs between $\p H_{2r'}$ and $\p H_{r''}$,
but now we require that the interfaces be well
separated on the inner boundary $\p H_{2r'}$,
instead of on the outer boundary.
As in the proof of Lemma~\ref{l.improve},
it is not too hard to see that conditioned on these
events there is probability bounded away from zero
(by a function of $j$)
that the crossings between $\p H_r$ and $\p H_{r'/2}$
will hook up nicely with the crossings between
$\p H_{2r'}$ and $\p H_{r''}$.
Basically, we only need to arrange that the channels
$\hat\beta_i$ for the inner crossings will cross the corresponding
channels of the outer crossings.  
The proof of the right hand inequality in~\eref{e.qm}
now follows from Lemma~\ref{l.plenty} and the corresponding
statement for crossings with interfaces well-separated in the
inner boundary, which is proved in the same way.
If $r''\le8\,r'$, then the right hand inequality
in~\eref{e.qm} is a consequence of Lemmas~\ref{l.improve}
and~\ref{l.plenty}. A similar proof applies when
$r'\le 8\, r$.

It now remains to prove the left hand inequality in~\eref{e.qm}.
If $r''<2\,r'$, then $\Pb{A_j(r',r'')}$ is bounded away from
zero (if we assume $\Pb{A_j(r,r'')}>0$)
and we are done since $\Pb{A_j(r,r'')}\le \Pb{A_j(r,r')}$.
Otherwise, we argue that
$\Pb{A_j(r,r'')}\le \Pb{A_j(r,r')}\,\Pb{A_j(2\,r',r'')}$,
by independence on disjoint subsets,
and $\Pb{A_j(r',2\,r')}\,\Pb{A_j(2\,r',r'')}\le C\,\Pb{A_j(r',r'')}$,
by the right hand inequality in~\eref{e.qm}.
Since $\Pb{A_j(r',2\,r')}$ is bounded away from zero, the left hand
inequality now follows.
\QED

We now generalize Proposition~\ref{p.qm} to arbitrary sequences
of crossings.

\begin{proposition}\label{p.gencol}
Let $j\ge 1$ be an integer, and fix a color sequence $X\in\{\mathrm{black},\mathrm{white}\}^j$.
The probabilities for the existence of $j$ disjoint crossings
whose colors match this sequence in counterclockwise order
also satisfy the inequalities in Proposition~\ref{p.qm}.
The corresponding statement also holds
in the setting of critical bond percolation on the square grid.
\end{proposition}

\proof
We start with the easier case where all the colors in the sequence $X$
are the same, say black.
Suppose that $r'>2\,r$ and $r''>2\,r'$.
Consider the event that (a) there are $j$ disjoint black crossings
from $\p H_r$ to $\p H_{r'}$ and (b) there are $j$ disjoint
black crossings from $\p H_{r'}$ to $\p H_{r''}$
and (c) there are $j$ disjoint black circuits separating
$\p H_{r'/2}$ from $\p H_{r'}$
and (d) there are $j$ disjoint black circuits
separating $\p H_{r'}$ from $\p H_{2\,r'}$
and (e) there are $j$ disjoint crossings from
$\p H_{r'/2}$ to $\p H_{2\,r'}$.
Note that if we choose any one path in each of (a)--(e),
we can extract from the union a crossing from $\p H_r$ to 
$\p H_{r''}$. To see that we actually have
$j$ disjoint crossings from $\p H_r$ to $\p H_{r''}$, note that if
we remove any $j-1$ hexagons, then there is still one path
remaining in each of (a)--(e), and consequently, there is
still a crossing from $\p H_r$ to $\p H_{r''}$.
Thus, Menger's theorem
(see~\cite{Diestel})
implies that when (a)--(e) all hold
there are $j$ disjoint crossings from $\p H_r$ to 
$\p H_{r''}$. By the Harris-FKG inequality, the events (a)--(e)
are all positively correlated. By RSW, events
(c)--(e) have probabilities bounded away from zero
(assuming that (a) has positive probability).
The inequality corresponding to the right hand inequality in~\eref{e.qm}
now follows. The corresponding left hand inequality,
as well as the cases where $r'\le 2\,r$
or $r''\le 2\,r'$ are now proved as in the proof of Proposition~\ref{p.qm}.

We now assume that both colors appear in $X$,
and indicate the adaptations necessary in the proof of
Proposition~\ref{p.qm} to generalize to the present setting.
The quantity $s(r,R)$ needs to be defined slightly differently.
In the modified definition of $s(r,R)$, still only 
interfaces between crossings of opposite colors are considered. 
Suppose that $\gamma_1$ and $\gamma_2$ are two adjacent interfaces
from $\p H_r$ to $\p H_R$, and that $Q$ is the component
of $H_R\setminus (H_r\cup\gamma_1\cup\gamma_2)$ between them.
Let $\dist(z,Z;Q)$ denote the infimal length of a path from $z$ to $Z$
in $\closure Q$.
For $s>0$ set $W(s):=\{z\in Q:\dist(z,\p H_R;Q)\le s\}$.
The {\bf margin} between $\gamma_1$ and $\gamma_2$ is
defined as the supremum of the set of 
$s>0$ such that any path connecting
$\gamma_1$ and $\gamma_2$
in $W(s)$ has length at least $s$.
Now $s(r,R)$ is redefined as the smallest margin
among any two consecutive interfaces.

Lemma~\ref{l.apart} is still valid with this new definition of $s(r,R)$.
In fact, the only change needed in the proof is that instead of
looking for a white crossing 
in $Y_i\setminus B(z_i,\delta\,R)$ from
$\p\hat\beta_i$ to $\p H_R$, one looks for
a crossing in $Y_i\setminus B(z_i',2\,\delta\,R)$,
where $z_i'$ is the last point on the arc
$\p Y\cap\p H_R$, directed away from $\alpha_1$,
that has distance at most $\delta\,R$
from $\beta_i$.
(Here, we assume that $\delta<1/100$, say.)

We now explain how this new definition facilitates the
obvious analogue of Lemma~\ref{l.improve}.
Indeed, suppose that $\gamma_1$ and $\gamma_2$
are two adjacent interfaces, there are
at least $j_1\le j$ disjoint
black crossings in the sector
$Q$ of $H_R\setminus H_r$ between $\gamma_1$ and $\gamma_2$,
and the margin between $\gamma_1$ and $\gamma_2$ is at least $s$.
Let 
\begin{equation*}\begin{aligned}
Q_i& :=\{z\in W(s): 2\,i\,s/(2\,j_1)\le \dist(z,\p H_R;Q)\le (2\,i+1)\,s/(2\,j_1)\}
\,,\\
Q_i^* &:=\{z\in W(s): 2\,i\,s/(2\,j_1)\le \dist(z,\gamma_1;Q)\le (2\,i+1)\,s/(2\,j_1)\}
\,,
\end{aligned}\end{equation*}
where $W(s)$ is as above.
We may then consider the event that in each $Q_i$, $i=0,1,\dots,j_1-1$,
there is a black crossing from $\gamma_1$ to $\gamma_2$,
and in each $Q_i^*$, $i=0,1,\dots,j_1-1$
there is a black crossing from $\p H_R$ to 
$\p W(s)\setminus (\gamma_1\cup\gamma_2\cup \p H_R)$,
and moreover, the latter crossings continue through well directed
channels all the way to $\p H_{4R}$, as in the proof of
Lemma~\ref{l.improve}.
An application of Menger's theorem,
as in the monochromatic case above, will then guarantee
that at the end $j_1$ disjoint black crossings between $\gamma_1$ and
$\gamma_2$ will extend all the way to 
$H_{4R}$. A compatible construction is applied to
each of the other pairs of adjacent interfaces.

Of course, when we condition on the interfaces
$\gamma_0,\gamma_1,\dots,\gamma_{k-1}$, we do not
know how many crossings we will have between each pair of
adjacent interfaces.
 But $k=O(R/s(r,R))$, and so the number of distinct
patterns in which crossings with color sequence type $X$
can occur is bounded by a function of $\delta$ and $j$.
(Specifically, a pattern for $X$ is a specification of
how many crossings are selected between each pair of adjacent interfaces
to make up the sequence of crossings compatible with $X$.
There may very well be additional crossings that are ignored.)
Thus, the most likely pattern given the interfaces
occurs with probability bounded away from zero given that
there are crossings of color sequence $X$
and the construction may be based on this most likely pattern.
The occurrence of this pattern given the interfaces and
the information about the color of hexagons adjacent to
the right hand sides of the interfaces will be a monotone
function in the collection of white hexagons in the regions
between interfaces that have white hexagons on their boundary,
and monotone in black hexagons in the other regions.
Thus, again, the Harris-FKG inequality is applicable.
(We do not want to condition on the exact number
of crossings between two adjacent interfaces, as this is not
a monotone function of the configuration.)

Similarly, when we attempt to glue crossings between two
different annuli, we condition on the interfaces, and then
aim for the most likely pattern in each annulus.
These are essentially the only modifications needed in the proof.
\QED

\begin{remark}\label{r.qmc}
The analogous statements for critical percolation in cones and
wedges also holds, with similar proofs. The wedge case is, in
fact, easier.
\end{remark}

\begin{remark}\label{r.primus}
It is also clear that the above proof shows that if we prescribe
$j$ specific disjoint arcs on $\p H_R$ and require the crossings
from $\p H_r$ to land on
these arcs, with a prescribed color for every arc,
the probability for this event is at least a positive constant
times the probability to have $j$ crossings with this
sequential color pattern, where the constant only depends
on the smallest angle at $0$ subtended by any of the $j$ arcs
(provided that $R$ is not too small, so that each of the arcs
has at least one hexagon unshared with any other arc, say).
\end{remark}

As a further application, we prove the following
result about $5$-arm and $6$-arm exponents in $\Z^2$.

\begin{corollary}\label{c.56}
Consider critical bond percolation on $\Z^2$.
For $R>r\ge 1$ let $F(r,R)$ denote the event that 
there are five open crossings between distances $r$ and $R$
from $0$ of types primal, primal, dual, primal, dual,
in circular order.
Then
\begin{equation}\label{e.fivearm}
 C^{-1}(r/R)^{2}\le
\Pb{F(r,R)} \le C\,(r/R)^{2},
\end{equation}
where $C>0$ is a universal constant.
Moreover, the probability that there are $6$ alternating
crossings: primal, dual, primal, dual, primal, dual,
between distances $r$ and $R$ is at most
$C\,(r/R)^{2+\eps}$, for some constant $\eps>0$.
The same statement applies to any sequence obtained
by inserting one additional primal or dual
entry to the list (primal, primal, dual, primal, dual).
\end{corollary}

This result is essentially due to~\cite[Lemma 5]{KSZ} (in
the context of site percolation on the triangular grid,
though the proof is equally applicable to $\Z^2$).
They omit some of the details,
because the proof is long and the argument is similar to
the proof of~\cite[Lemma 4]{K2}.
Now, we can easily present an essentially complete
and relatively short argument.

\proof
We begin with the basic argument from~\cite{KSZ}.
Divide the boundary of the circle $\p B(0,R)$ into
$5$ equal arcs, $A_1,\dots,A_5$, in counterclockwise order.
For concreteness, let's take each $A_j$ as the arc between
angles $(2\,j-1)\,\pi/5$ and $(2\,j+1)\,\pi/5$.
For a vertex $v\in B(0,R)$,
let $F_v$ be the event that
there are primal (open) crossings from $v$ to $A_1,A_3$ and $A_4$ and
dual crossings from dual vertices adjacent to $v$ to $A_2$ and to $A_5$, 
and the primal crossings are disjoint, except at $v$.
(By planarity, it follows that the dual crossings are disjoint.)
We claim that $F_v$ can happen for at most one vertex in
$B(0,R/2)$. Indeed, suppose that $F_v\cap F_u$ holds,
where $v,u$ are vertices in $B(0,R/2)$.
Let $\alpha_i$, $i=1,3,4$ denote some simple primal crossings from $v$
to the arcs $A_i$, which are disjoint, except at $v$.
Similarly, let $\alpha_i'$, $i=1,3,4$, be the corresponding
paths for $u$.
Since $\alpha_1\cup\alpha_3$ separates
$A_2$ from $A_5$ in $B(0,R)$, it is clear that
$u\in \alpha_1\cup\alpha_3$. Similarly,
$u\in\alpha_1\cup\alpha_4$. Since $\alpha_3\cap\alpha_4=\{v\}\subset\alpha_1$,
it follows that $u\in \alpha_1$. Let $\beta_1$ be the arc of $\alpha_1$
from $u$ to $A_1$. Since  $\beta_1\cup\alpha_3'$ is a path from $A_1$ to $A_3$,
it contains $v$ or separates $v$ from $A_2$ or from $A_5$.
The latter two possibilities are ruled out by the dual crossings
to $A_2$ and $A_5$ starting at dual vertices adjacent to $v$.
Thus, $v\in\beta_1\cup\alpha_3'$, and similarly,
$v\in \beta_1\cup\alpha_4'$. But since $\alpha_3'\cap\alpha_4'=\{u\}\subset\beta_1$,
we conclude that $v\in \beta_1$, which implies $u=v$.

We now claim that $F:=\bigcup_{v\in B(0,R/2)} F_v$ has probability bounded away from $0$. 
Consider the event that there is a crossing from $A_1$ to $A_4$, and consider the rightmost
such crossing $\ell$ (in the sense that it separates any other crossing from $A_5$).
If there is an open path from $A_3$ to $\ell$, but there is no open path from $A_3$
to $A_1$ disjoint from $\ell$, then $F_w$ will hold, where
$w$ is the first vertex $v$ along $\ell$ (when $\ell$ goes from $A_1$ to $A_4$)
that connects to $A_3$ in the complement of $\ell$.
Thus, we need to show that there is probability bounded away from zero that this
happens with $w\in B(0,R/2)$.
Let $L_1$ be the line passing through the origin and the midpoint of $A_1$.
Let $L_2$ and $L_3$ be lines parallel with $L_1$ at distance $R/20$ and $R/10$
from $L_1$, on the side of $L_1$ that contains $A_5$.
By RSW, there is probability bounded away from zero for the existence of
a dual-open crossing from $A_1$ to $A_4$ in the strip between $L_2$ and $L_3$.
By conditioning on the leftmost such crossing (the one closest to $L_1$),
it is easy to see that there is probability bounded away from
zero that such a dual crossing exists and is also connected to $A_5$ by a
dual-open path.
If moreover we have a primal crossing from $A_1$ to $A_4$ in the strip
between $L_1$ and $L_2$, which happens with
probability bounded away from zero, then the rightmost primal crossing between 
$A_1$ and $A_4$ will be contained in the strip
between $L_1$ and $L_3$.
Conditioned on the latter event and on the latter crossing $\ell$,
it happens with probability bounded away from zero that
there is a primal crossing from $A_3$ to $\ell$ whose endpoint on $\ell$
(which is its only point on $\ell$) is within distance
$R/5$ of the origin and there is a dual crossing from $A_2$ to a dual
vertex adjacent to $\ell$ that is within distance $R/5$ of the origin.
On that event, $F$ holds. Thus, $\Pb{F}$ is bounded away from zero.

It is easy to see that the proof
of Proposition~\ref{p.gencol} implies
that $\Pb{F_v}\le O(1)\,\Pb{F_w}$ for $v,w\in B(0,R/2)$.
Since the events $F_v$ are disjoint, and since their sum is of order $1$,
it follows that each $F_v$, $v\in B(0,R/2)$ has probability of order $R^{-2}$.
In particular $R^2\,\Pb{F_0}$ is bounded away from zero and $\infty$.
Now~\eref{e.fivearm} follows from Remark~\ref{r.primus}.

The statements regarding the $6$-arm exponent now follow
from Reimer's inequality~\cite{Reimer}. Alternatively, we may also deduce them
from Remark~\ref{r.primus}, as follows. If we fix arcs 
$A_1,\dots,A_6$ in counterclockwise order on the radius $R$ circle, 
where the crossings are required to land,
and we require a primal crossing to $A_1$ and a dual crossing to $A_2$,
then we may condition on the most counterclockwise primal crossing $\gamma$ connecting to $A_1$,
then sequentially on the most clockwise crossings to $A_2,A_3,\dots,A_5$  of the required type.
The conditional probability for having yet another crossing to $A_6$ is 
still bounded by $O(1)\,(r/R)^\eps$, for
some constant $\eps>0$. Now we may apply Remark~\ref{r.primus}, to complete the proof.
\QED

\bigskip
\noindent{\bf Acknowledgments}.
We thank Harry Kesten for useful advice.

\bigskip
\filbreak
\begingroup
{
\small
\parindent=0pt

Microsoft Corporation\\
One Microsoft Way\\
Redmond, WA 98052, USA\\ 
{\tt
{http://research.microsoft.com/\string~schramm/}}

\bigskip

Department of Mathematics \\
Chalmers University of Technology 
and G\"{o}teborg University \\
Gothenburg, Sweden \\
{\tt
{steif@math.chalmers.se} \\
{http://www.math.chalmers.se/$\sim$steif/}}
}

\filbreak

\endgroup

\end{document}